\documentclass[preprint,12pt,number]{elsarticle}

\usepackage[paper=letterpaper, margin=1in]{geometry}
\usepackage{amssymb}
\usepackage{amsmath}
\usepackage{graphicx}
\usepackage{booktabs}
\usepackage{hyperref}

\hypersetup{
  pdftitle={Finite Element Eigenfunction Network (FEENet): A Hybrid Framework for Solving PDEs on Complex Geometries},
  pdfauthor={Shiyuan Li; Hossein Salahshoor},
  pdfkeywords={FEENet, Neural operator, Finite element method, Eigenfunction, DeepONet}
}

\usepackage{caption}
\usepackage{multirow}
\usepackage{makecell}
\usepackage{tabularx}
\usepackage{subcaption}
\usepackage{amsthm}
\newtheorem{theorem}{Theorem}

\journal{}

\begin{document}

\makeatletter
\def\ps@pprintTitle{%
  \let\@oddhead\@empty
  \let\@evenhead\@empty
  \let\@oddfoot\@empty
  \let\@evenfoot\@empty}
\makeatother

\begin{frontmatter}

\title{Finite Element Eigenfunction Network (FEENet): A Hybrid Framework for Solving PDEs on Complex Geometries}

\author[cee]{Shiyuan Li}
\author[cee,mems]{Hossein Salahshoor\corref{cor1}}
\cortext[cor1]{Corresponding author}
\ead{hossein.salahshoor@duke.edu}

\affiliation[cee]{organization={Department of Civil and Environmental Engineering},
    addressline={Duke University},
    city={Durham},
    state={NC},
    country={USA}}

\affiliation[mems]{organization={Department of Mechanical Engineering and Materials Science},
    addressline={Duke University},
    city={Durham},
    state={NC},
    country={USA}}

\begin{abstract} 
Neural operators aim to learn mappings between infinite-dimensional function spaces, but their performance often degrades on complex or irregular geometries due to the lack of geometry-aware representations. We propose the \textit{Finite Element Eigenfunction Network} (FEENet), a hybrid spectral learning framework grounded in the eigenfunction theory of differential operators. For a given domain, FEENet leverages the Finite Element Method (FEM) to perform a one-time computation of an eigenfunction basis intrinsic to the geometry. PDE solutions are subsequently represented in this geometry-adapted basis, and learning is reduced to predicting the corresponding spectral coefficients. Numerical experiments conducted across a range of parameterized PDEs and complex two- and three-dimensional geometries, including benchmarks against the seminal DeepONet framework \cite{lu2019deeponet}, demonstrate that FEENet consistently achieves superior accuracy and computational efficiency. We further highlight key advantages of the proposed approach, including resolution-independent inference, interpretability, and natural generalization to nonlocal operators defined as functions of differential operators. We envision that hybrid approaches of this form, which combine structure-preserving numerical methods with data-driven learning, offer a promising pathway toward solving real-world PDE problems on complex geometries.
\end{abstract}


\end{frontmatter}

\section{Introduction}

Partial differential equations (PDEs) serve as the mathematical foundation of a wide range of models in science and engineering, governing phenomena such as elasticity, fluid flow, heat transfer, electromagnetism, and wave propagation \cite{evans2022partial, brezis2011functional}. In practical applications, PDEs must be solved on domains characterized by complex and irregular geometries with intricate structures \cite{li2009solving}, with examples including porous and composite materials \cite{torquato1991random, milton2022theory}, biological tissues \cite{peskin1977numerical, humphrey2004introduction, xu2020phase, henninger2010validation, salahshoor2020transcranial, salahshoor2022mechanics}, architected and additively manufactured structures \cite{hussein2014dynamics, deshpande2001effective,ulloa2025homogenized, salahshoor2018material,  mahabadi2026data}, and geophysics \cite{virieux2011review, li2009solving, assimaki2003effects, veveakis2007thermoporomechanics}. In such settings, the geometry of the domain plays a central role in shaping the PDE solutions, and, consequently, any computational framework for solving PDEs in realistic scenarios must be able to faithfully account for domain complexity while maintaining reasonable accuracy and stability.

For several decades, the FEM has served as the gold standard for the numerical solution of PDEs, particularly in complex geometries, where it provides strong theoretical guarantees in terms of stability, convergence, and approximation accuracy \cite{zienkiewicz1977finite,hughes2003finite, ciarlet2002finite}. To alleviate the challenges posed by complex or evolving domains, a broad class of numerical techniques has been developed, including fictitious domain methods \cite{glowinski1994fictitious}, composite and extended finite element methods \cite{moes1999finite, melenk1996partition, dolbow2009efficient}, isogeometric analysis \cite{hughes2005isogeometric, bazilevs2023computational}, immersed and embedded interface methods \cite{leveque1994immersed,li2006immersed, zhang2004immersed}, embedded boundary and cut-cell methods \cite{peskin2002immersed, colella1998cartesian, johansen2008embedded, main2018shifted, atallah2022high}, as well as diffuse or phase-field domain approaches in which geometry is represented implicitly through auxiliary fields \cite{li2009diffuse}. While these methods substantially expand the range of geometries that can be treated within the FEM framework, they typically introduce additional stabilization mechanisms, interface treatments, or localized refinement near boundaries, leading to increased algorithmic and computational complexity. Consequently, despite the maturity and robustness of FEM-based approaches, the reliable treatment of highly intricate or evolving geometries is difficult and costly, especially in three dimensions, and high-speed calculations demanded in practice remain a significant computational challenge.

In many contemporary applications in design and engineering, PDEs must be solved repeatedly for varying inputs, parameters, or forcing terms. In such many-query regimes, the repeated execution of high-fidelity numerical solvers such as FEM can become computationally prohibitive, motivating the development of surrogate models that approximate the underlying solution operator. Neural operators offer an alternative paradigm by directly learning mappings between infinite-dimensional function spaces, thereby approximating the PDE solution operator itself rather than individual solutions \cite{khoo2021solving}. Building on early universal approximation results for operators \cite{chen1993approximations,chen1995universal}, the introduction of DeepONet \cite{lu2019deeponet} established a practical framework for operator learning and spurred rapid advances in the field. Subsequent developments, including variants and extensions of DeepONet \cite{jin2022mionet,zhang2023belnet,lu2022comprehensive,zhang2025discretization,zhu2023reliable,howard2022multifidelity,lin2022learning,lin2023b,zhang2024d2no} and Fourier Neural Operators (FNOs) and their extensions \cite{li2020fourier,li2020neural,li2024physics,pathak2022fourcastnet,li2023solving,wen2022u,li2023fourier}, have demonstrated impressive performance across a wide range of scientific and engineering problems \cite{zhang2025ap, wen2022u, lehmann20243d}. Despite these successes, most existing neural operator architectures rely on fixed, geometry-agnostic representations—such as uniform grids, Fourier bases, or global coordinate embeddings—which limits their ability to generalize to domains with complex or irregular geometries, precisely the regimes where FEM remains most effective.

To address the geometric limitations of standard neural operator architectures, several recent works have proposed extensions tailored to complex domains. For instance, \cite{he2024geom} introduced a geometry-aware adaptation of DeepONet for three-dimensional shapes. A related line of research seeks to decouple the representation of the solution space from operator learning by explicitly learning or prescribing basis functions; see, e.g., \cite{ingebrand2025basis,lee2024training}. Along similar lines, POD-DeepONet \cite{lu2022comprehensive} employs proper orthogonal decomposition of training data to construct a reduced basis, while the fixed-basis coefficient-to-coefficient operator network (FB-C2CNet) \cite{chen2025learning} assumes a prescribed basis and learns the corresponding coefficient mappings. While these approaches represent important progress toward geometry-aware operator learning, the resulting bases are typically data-dependent or weakly tied to the underlying differential operator, which can limit robustness and generalization on complex geometries.

In this work, we propose a universal hybrid framework, the \textit{Finite Element Eigenfunction Network} (FEENet), for solving PDEs on complex geometries. FEENet combines the geometric robustness and spectral structure of FEM with the efficiency of neural operator learning. Building on the two-network architecture of DeepONet \cite{lu2019deeponet}, we replace the learned trunk representation with an operator-adapted basis obtained by solving the eigenvalue problem associated with the governing differential operator via FEM. By exploiting the completeness of elliptic eigenfunctions, solutions are represented in this precomputed basis, and learning is reduced to predicting the corresponding spectral coefficients. This separation avoids the insufficient information fusion that can arise when basis functions and coefficients are learned simultaneously \cite{wang2022improved}. Conceptually, the finite element eigensolver acts as an \emph{a priori} trunk network equipped with an optimal, geometry-aware representation, which explains the rapid convergence and stability of the proposed method. For a fixed geometry, FEENet then requires only a one-time eigenvalue computation, after which the resulting eigenfunctions provide a robust basis for the solution spaces of a broad class of steady-state and time-dependent PDEs. The framework preserves resolution-independent inference, supports queries at arbitrary spatial locations, and significantly reduces training cost, with accuracy governed primarily by the fidelity of the finite element spectral approximation rather than by ad hoc network architecture or hyperparameter tuning.

The remainder of this paper is organized as follows. Section~\ref{sec:FEENet} presents the mathematical formulation and architecture of FEENet: After reviewing a key result from functional analysis concerning eigenfunction expansions of differential operators, we introduce the core ideas underlying the proposed framework and describe its implementation. Section~\ref{sec:experiments} reports numerical experiments on Poisson and heat problems posed on three representative geometries of increasing complexity, and compares the performance of FEENet with DeepONet (MIONet). Section~\ref{sec:discussion} provides a detailed discussion of the results, highlighting salient features of the method such as resolution-independent inference and computational efficiency, and outlines potential extensions and directions for future work.

\section{Finite Element Eigenfunction Network}
\label{sec:FEENet}
This section delineates the central idea, mathematical foundations, and architectural design of the proposed Finite Element Eigenfunction Network (FEENet) method for approximating solutions of partial differential equations (PDEs). 

\subsection{Preliminaries: Spectral representation of elliptic operators}

Let $\Omega \subset \mathbb{R}^d$ be a bounded domain with a sufficiently smooth boundary, and let $\mathcal{L}_\mu$ be a linear, self-adjoint, strongly elliptic differential operator of order $2m$ on $\Omega$, endowed with homogeneous Dirichlet boundary conditions. We consider elliptic problems of the form
\begin{equation}
\mathcal{L}_\mu u = f \quad \text{in } \Omega,
\qquad
u = 0 \quad \text{on } \partial\Omega,
\label{eq:elliptic_pde}
\end{equation}
where $\mu$ denotes possible PDE parameters, $f$ is a given forcing term, and $u$ is the solution field. Our objective is to characterize the associated solution operator $\mathcal{G}: f \mapsto u$. In the following spectral discussion, we omit the parameter $\mu$ from the notation for simplicity.

An eigenfunction of $\mathcal{L}$ is a nontrivial solution $\psi$ of the spectral problem
\begin{equation}
\mathcal{L}\psi = \lambda \psi \quad \text{in } \Omega,
\qquad
\psi = 0 \quad \text{on } \partial\Omega,
\label{eq:eigen_problem}
\end{equation}
where $\lambda \in \mathbb{R}$ is the associated eigenvalue. Classical elliptic theory ensures that the spectrum of $\mathcal{L}$ is discrete, real, and unbounded above, and that the corresponding eigenfunctions are smooth in $\Omega$.

A fundamental result in the spectral theory of elliptic operators states that the eigenfunctions of $\mathcal{L}$ form a complete system not only in $L^2(\Omega)$, but also in the natural energy space $H^m(\Omega)$, the Sobolev space of functions that are weakly differentiable up to order $m$. This result was first established by Browder for self-adjoint elliptic operators of arbitrary order, and was subsequently improved by G{\aa}rding and Agmon \cite{Garding1953,Agmon1962}.

\begin{theorem}[Completeness of elliptic eigenfunctions {\cite{Browder1952,Browder1953}}]
Let $\mathcal{L}$ be a self-adjoint linear elliptic differential operator of order $2m$ on a bounded domain $\Omega \subset \mathbb{R}^d$ with homogeneous Dirichlet boundary conditions. Then the (real) eigenfunctions of $\mathcal{L}$ form a complete system in $H^m(\Omega)$. In particular, for any $u \in H^m(\Omega)$ and any $\varepsilon>0$, there exists a finite linear combination of eigenfunctions that approximates $u$ in the $H^m$-norm.
\end{theorem}

As a consequence, any solution $u \in H^m(\Omega)$ of \eqref{eq:elliptic_pde} admits a representation of the form
\begin{equation}
u = \sum_{j=1}^{\infty} c_j \phi_j,
\qquad
c_j = \langle u,\phi_j\rangle_{L^2(\Omega)},
\label{eq:SchauderRep}
\end{equation}
where $\{\phi_j\}_{j\ge1}$ are eigenfunctions of $\mathcal{L}$, and the series converges in the $H^m(\Omega)$ norm.

\subsection{Illustration of the FEENet idea}

\begin{figure}[htbp]
    \centering
    \includegraphics[width=\linewidth]{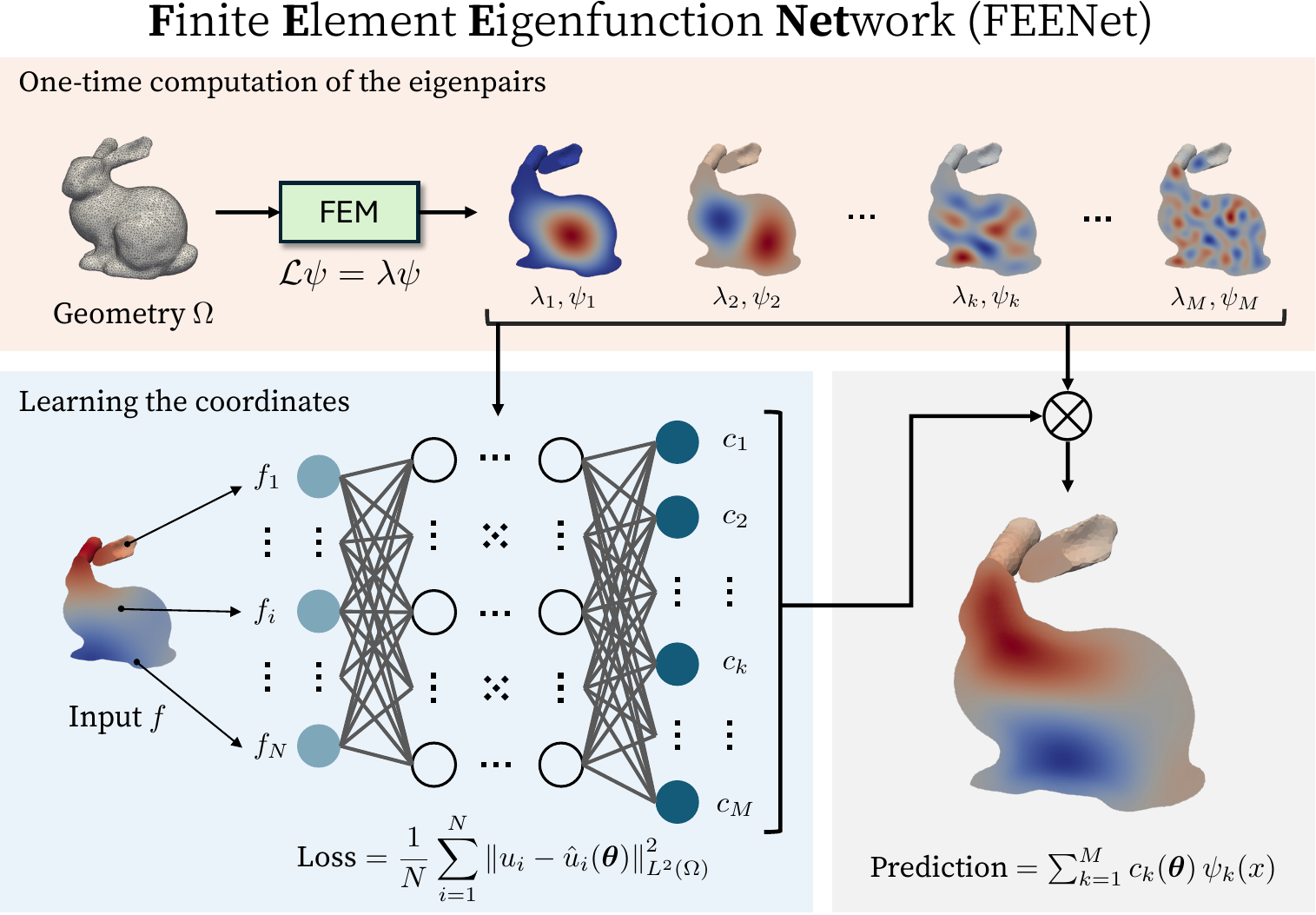}
    \caption{\textbf{Finite Element Eigenfunction Network (FEENet).}
    This framework learns the solution operator of a differential equation governed by
    operator $\mathcal{L}$ by decoupling the problem into two components:
    (\textit{Top}) An offline FEM solver pre-computes the eigenvalue--eigenfunction
    pairs $\{\lambda_k, \psi_k\}_{k=1}^M$ by solving the eigenvalue problem
    $\mathcal{L}\psi_k = \lambda_k \psi_k$ in $\Omega$ with homogeneous boundary conditions.
    These eigenfunctions serve as fixed, geometry-aware basis functions.
    (\textit{Bottom}) For each training sample, a branch network maps the input
    function $f$ to spectral coordinates $\{c_k(\boldsymbol{\theta})\}_{k=1}^M$,
    and the predicted solution is reconstructed via a linear combination of basis
    functions and corresponding learned coordinates.}
    \label{fig:hero}
\end{figure}

Taken together, these results justify the representation of solutions to elliptic partial differential equations on arbitrary bounded domains as convergent expansions in eigenfunctions of the governing operator. This spectral viewpoint provides a geometry-aware and operator-adapted coordinate system for solution spaces. 
The central idea of our framework is to exploit this spectral representation by precomputing finite-element approximations of the eigenfunctions of $\mathcal{L}$ and learning solution coordinates $\{c_j\}$ in this operator-adapted basis via a neural network. In this manner, FEM yields the geometry-aware spectral basis, and the neural network learns the solution operator in those coordinates. Figure \ref{fig:hero} schematically summarizes the idea.

Let us utilize the Laplace-Beltrami operator as a canonical example to highlight the distinct features of FEENet. For a given geometry, one can use the Laplacian eigenfunctions for all differential operators. Of course, if $\mathcal{L}$ is a differential operator of order higher than 2, using eigenfunctions corresponding to $\mathcal{L}$ would be more efficient.

In FEENet, we leverage the finite element method (FEM) to pre-compute the eigenfunctions of the Laplacian, solving the eigenvalue problem:
\begin{equation}
    \begin{cases}
        -\Delta \phi_k = \lambda_k \phi_k & \text{in } \Omega, \\
        \mathcal{B}(\phi_k) = 0 & \text{on } \partial\Omega,
    \end{cases}
\end{equation}

where $\mathcal{B}$ denotes appropriate boundary conditions. These eigenfunctions $\Phi = \{\phi_k\}_{k=1}^M$ form a robust, fixed orthogonal basis for $H^1(\Omega)$.

As demonstrated by \cite{reuter2006laplace}, the spectrum $\{\lambda_k\}$, sometimes referred to shape-DNA, uniquely characterizes the geometry. By employing these geometry-aware modes as the trunk functions, we ensure that the predicted solution $u(x)$ naturally respects the domain's structure. In our proposed hybrid framework, learning the coordinates of the solution is transferred to the neural network (Branch Net), which learns the mapping (forcing, boundary conditions, parameters) $\mapsto \{c_k\}_{k=1}^M$. This design decouples the encoding of geometric complexity (via the offline Laplacian solver) from the learning of physical laws (via the neural network).

\subsection{Remark on the connection of FEENet with existing operator learning frameworks}
FEENet, in a sense that is elaborated in this subsection, builds upon the theory of DeepONet \citep{lu2019deeponet}. In the standard DeepONet framework \citep{lu2019deeponet}, the solution operator $\mathcal{G}: f \mapsto u$ is approximated by a sum of products of two neural networks:
\begin{equation}
    u(x) \approx \mathcal{G}(f)(x) = \sum_{k=1}^p \underbrace{b_k(f(x_1), \dots, f(x_m))}_{\text{Branch Net}} \cdot \underbrace{t_k(x)}_{\text{Trunk Net}}, \quad x \in \Omega.
\end{equation}
The \textit{branch network} encodes the discrete input function $f$ evaluated at fixed sensor locations $\{x_i\}_{i=1}^m$, while the \textit{trunk network} learns a set of basis functions $\{t_k(x)\}$ over the domain $\Omega$. While this architecture is flexible, the data-driven trunk network does not incorporate explicit geometric priors. Consequently, relying purely on coordinate-based learning can be inefficient, often requiring substantial training data to accurately resolve high-frequency features or satisfy complex boundary conditions \citep{kontolati2023influence}.

\subsection{FEENet Architecture}
\label{sec:architecture}

FEENet replaces the learned trunk network with a fixed projection onto the first $M$ eigenfunctions. For a given geometry, the architecture consists of:
\begin{enumerate}
    \item \textbf{Preprocessing:} A one-time computation of eigenpairs $(\lambda_k, \phi_k)$ using a FEM solver (e.g., DOLFINx \cite{baratta2023dolfinx}).
    \item \textbf{Branch Network:} A lightweight neural network that maps the discretized input $f(x_i)$ to a set of spectral coordinates $\mathbf{c} = [c_1, \dots, c_M]^T$.
    \item \textbf{Spectral Synthesis:} The final solution is reconstructed as $u(x) = \sum_{k=1}^M c_k \hat{\phi}_k(x)$, where $\hat{\phi}_k$ may represent normalized eigenfunctions.
\end{enumerate}

While we have developed this framework independently, the authors remark that the formulation is similar to the recent FB-C2CNet paradigm \citep{chen2025learning}, but specifically tailors the basis to the eigenfunction along with the physical domain through the finite element method. For time-dependent problems, we further incorporate the analytical spectral decay $e^{-\lambda_k t}$ into the computation graph, as detailed in Section~\ref{sec:experiments}.

\subsection{Training Objective and Loss Function}

The prediction of FEENet is reconstructed using a truncated spectral expansion based on the eigenfunctions of the differential operator $\mathcal{L}$ defined in Equation~\eqref{eq:eigen_problem}.

For each training sample indexed by $i$, the predicted solution is expressed as:
\begin{equation}
\label{eq:reconstruction_general}
\hat{u}_i(x; \boldsymbol{\theta}) = \sum_{k=1}^M c_{k,i}(\boldsymbol{\theta})\, \psi_k(x),
\end{equation}
where $\{\lambda_k, \psi_k\}_{k=1}^M$ denote the first $M$ eigenvalue-eigenfunction pairs of $\mathcal{L}$, and $c_{k,i}(\boldsymbol{\theta})$ are the spectral coordinates predicted by the branch network.

The parameters $\boldsymbol{\theta}$ of the branch network are optimized by minimizing the mean-squared error loss:
\begin{equation}
\label{eq:loss}
\min_{\boldsymbol{\theta}} \mathcal{L}(\boldsymbol{\theta})
=
\frac{1}{N}
\sum_{i=1}^{N}
\left\|
u_i - \hat{u}_i(\boldsymbol{\theta})
\right\|_{L^2(\Omega)}^2,
\end{equation}
where $N$ denotes the number of training samples, $u_i$ is the ground-truth solution, and $\hat{u}_i(\boldsymbol{\theta})$ is the corresponding network prediction.

The loss function is minimized using the Adam optimizer \cite{kingma2017adammethodstochasticoptimization}. Importantly, the eigenfunctions $\{\psi_k\}$ and eigenvalues $\{\lambda_k\}$ are fixed throughout training and are not treated as trainable parameters. The branch network therefore only learns the mapping from input functions to spectral coordinates, which significantly reduces the optimization complexity compared to the standard DeepONet architectures.

For time-dependent problems or systems with analytical structure (e.g., forced dynamics), additional physical constraints can be incorporated directly into the reconstruction formula. The specific forms for different PDE classes are detailed in Section~\ref{sec:experiments}.

\section{Numerical Experiments}
\label{sec:experiments}
This section evaluates the proposed FEENet framework across three PDE benchmarks defined on diverse geometries and dimensions. The tested geometries include a 2D Unit Square, a 2D Fins, and the 3D Bunny \citep{StanfordBunny1994}. The benchmark problems encompass the Poisson problem, the homogeneous heat problem, and the inhomogeneous heat problem with a time-independent forcing term.

We employ the FEniCSx software stack—DOLFINx, Basix, and UFL \cite{baratta2023dolfinx, scroggs2022construction, scroggs2022basix, alnaes2014unified}—to compute all high-fidelity ground truth solutions. To obtain diverse and physically meaningful test cases, we sample input functions (forcing terms and initial conditions) from Gaussian Random Fields (GRFs) with a Gaussian covariance kernel using the GSTools library \citep{muller2022gstools}. The variance was fixed to $\sigma^2 = 15$, while the correlation length scale was adjusted for each geometry to ensure sufficient spatial variability: $\ell=0.3$ for the Unit Square, $\ell=0.15$ for the Fins, and $\ell=0.4$ for the Bunny. The dataset contains 2,000 samples for each PDE problem on each geometry. Detailed information on the GRF generation procedure, including the mathematical formulation and the Randomization method, is provided in~\ref{app:Gaussian_generation}.

\subsection{Training and Evaluation Protocols}
For model training, the DeepONet or its multi-input variant MIONet, and our proposed FEENet are implemented in the DeepXDE framework \cite{lu2021deepxde}. We denote network architectures using the notation $[N_{in}, h_1, \dots, h_k, N_{out}]$, 
where integers represent neuron counts in input, hidden, and output layers, respectively. For example, $[P, 400]$ indicates a shallow network with input dimension $P$ and output dimension 400. Unless explicitly indicated otherwise, the networks are initialized using the Glorot normal (Xavier) initializer \cite{glorot2010understanding} and trained using the Adam optimizer \cite{kingma2017adammethodstochasticoptimization} for 100{,}000 iterations with a batch size of 256. The loss function is the Mean Squared Error (MSE) evaluated at the discrete sensor locations. All training procedures were performed on a workstation equipped with an Intel(R) Core(TM) i9-14900 CPU, utilizing seven cores.

For performance evaluation reported in Table~\ref{tab:training}, we compute relative $L^2(\Omega)$ and $H^1(\Omega)$ errors on the computational mesh. Since model predictions are defined at discrete sensor locations, both predicted and reference solutions are first mapped onto the domain using a finite element representation. The detailed error evaluation procedure is described in~\ref{app:error_evaluation}.

\subsection{Eigenvalue Problem Solving}
\label{sec:laplacian}
To have a fair comparison of the FEENet computational cost with other methods,
we also include the time required to solve the associated eigenvalue problem.
We employ the eigenfunctions of the Laplacian--Beltrami operator as the canonical
spectral basis for diffusion-dominated PDEs. Specifically, for each geometry
$\Omega$, we compute the first $M=400$ eigenpairs $\{\lambda_k, \phi_k\}_{k=1}^M$
by solving the homogeneous Dirichlet eigenvalue problem.
The corresponding computational cost is reported in
Table~\ref{tab:computation_cost_for_eigenpairs}.

\begin{table}[ht]
\centering
\small
\caption{\textbf{Computational cost of basis generation on different geometries.} The table reports the wall-clock time (mean $\pm$ std, in minutes) required for the one-time computation of the first 400 Laplace eigenpairs using the FEM solver.}
\label{tab:computation_cost_for_eigenpairs}
\begin{tabular}{l c}
\toprule
\textbf{Geometry} & \textbf{Time (min)} \\
\midrule
Square & 0.120 $\pm$ 0.002 \\
Fins   & 1.073 $\pm$ 0.017 \\
Bunny  & 6.529 $\pm$ 0.077 \\
\bottomrule
\end{tabular}
\end{table}

\subsection{Poisson Problem}
\label{poisson_problem}
We first consider the steady-state Poisson equation with homogeneous Dirichlet boundary conditions:
\begin{equation}
    -\Delta u(x) = f(x), \quad x \in \Omega;
    \qquad u|_{\partial\Omega} = 0.
\end{equation}
We utilize the pre-computed Laplacian spectral basis $\{\phi_k\}_{k=1}^{400}$ described in Section~\ref{sec:laplacian}. To reflect the energy structure of the Poisson operator $(-\Delta)^{-1}$,
we scale the Laplacian eigenfunctions as
$\hat{\phi}_k = \lambda_k^{-1/2}\phi_k$,
which naturally emphasizes low-frequency modes. The network prediction for sample $i$ takes the form:
\begin{equation}
\label{eq:poisson_reconstruction}
\hat{u}_i(x; \boldsymbol{\theta}) = \sum_{k=1}^{400} c_{k,i}(\boldsymbol{\theta})\, \hat{\phi}_k(x),
\end{equation}
where $c_{k,i}$ are the spectral coordinates predicted by the branch network.

The forcing terms $f(x)$ were sampled from the GRF distribution described in~\ref{app:Gaussian_generation}. Ground truth solutions were computed using a standard linear FEM solver. All inputs and outputs are standardized using Z-score normalization.

For FEENet, the dataset consists of:
\begin{equation}
    \mathcal{T}_{\text{FEENet}} = \left\{ f_i(x), u_i(x) \right\}_{i=1}^{2000},
\end{equation}
where the normalized eigenfunctions $\{\hat{\phi}_k\}_{k=1}^{400}$ serve as fixed basis functions. The branch network maps the discretized forcing term to predicted coordinates with architecture $[P, 400]$ where $P$ is the number of sensor points (locations where the input function is evaluated). For the Unit Square, $P = 1225$; for the Fins, $P = 8889$; and for the Bunny, $P = 40565$.

For DeepONet, the dataset is structured as:
\begin{equation}
    \mathcal{T}_{\text{DeepONet}} = \left\{ \left( f_i(x), x \right), u_i(x) \right\}_{i=1}^{2000}.
\end{equation}
The network consists of a branch network $[P, 400, 400, 400]$ processing the forcing term and a trunk network $[d, 400, 400, 400]$ encoding spatial coordinates, where $d=2$ for 2D problems and $d=3$ for the Bunny. Both networks use ReLU activation and output 400-dimensional features that are combined via element-wise multiplication.

Learning rates are set to $4\times 10^{-5}$ (Unit Square, Fins) and $8\times 10^{-6}$ (Bunny) for FEENet, and $8\times 10^{-5}$ (Unit Square, Fins) and $4\times 10^{-5}$ (Bunny) for DeepONet.

\subsection{Homogeneous Heat Problem}
Next, we investigate the homogeneous heat equation:
\begin{equation}
    \partial_t u = D \Delta u, \quad x \in \Omega, \ t \in (0,T];
    \qquad u(x,0)=u_0(x).
\end{equation}

Using the same Laplacian spectral basis from Section~\ref{sec:laplacian}, the solution admits the analytical expansion $u(x,t) = \sum_{k=1}^M c_k(0) e^{-D\lambda_k t} \phi_k(x)$, where $c_k(0)$ are the initial spectral coordinates determined by the initial condition. Incorporating this temporal structure into the general reconstruction formula from Equation~\eqref{eq:reconstruction_general}, the network prediction becomes:
\begin{equation}
\label{eq:heat_homo_reconstruction}
\hat{u}_i(x,t; \boldsymbol{\theta}) = \sum_{k=1}^{400} c_{k,i}(0; \boldsymbol{\theta})\, e^{-D\lambda_k t}\, \phi_k(x),
\end{equation}
where the branch network learns the mapping from $u_0(x)$ to the initial coordinates $c_{k,i}(0)$, while the temporal decay $e^{-D\lambda_k t}$ is incorporated analytically.

The diffusion coefficient was set to $D=0.02$, and the system was simulated over the time interval $[0, T]$ with $T = 1.0$. Initial conditions were sampled from GRFs and post-processed to satisfy homogeneous Dirichlet boundary conditions; see~\ref{app:Gaussian_generation} for details. Ground truth trajectories were computed using the implicit Euler method with time step $dt=0.0025$. Snapshots were saved at 10 uniformly spaced time points from $t=0.1$ to $t=1.0$.

For FEENet, the dataset is structured as:
\begin{equation}
    \mathcal{T}_{\text{FEENet}} = \left\{ \left( u_{0,i}(x), t_j \right), u_i(x, t_j) \right\}_{i=1,j=1}^{2000,10},
\end{equation}
where the eigenpairs $\{\lambda_k, \phi_k\}_{k=1}^{400}$ are pre-computed and embedded in the model architecture. The branch network has architecture $[P, 400]$.

For DeepONet, the dataset takes the form:
\begin{equation}
    \mathcal{T}_{\text{DeepONet}} = \left\{ \left( u_{0,i}(x), (x,t) \right), u_i(x,t) \right\}_{i=1}^{2000}.
\end{equation}
The branch network processes initial conditions with architecture $[P, 400, 400, 400]$, while the trunk network encodes spatiotemporal coordinates $(x,t)$ with architecture $[d+1, 400, 400, 400]$.

Learning rates are $4\times 10^{-5}$ (Unit Square, Fins) and $2\times 10^{-6}$ (Bunny) for FEENet, and $5\times 10^{-5}$ (Unit Square, Fins) and $1\times 10^{-5}$ (Bunny) for DeepONet.

\subsection{Inhomogeneous Heat Problem}
Finally, we investigate the inhomogeneous heat equation with a time-independent forcing term:
\begin{equation}
    \partial_t u = D \Delta u + f(x), \quad x \in \Omega;
    \qquad u(x,0)=u_0(x).
\end{equation}

Projecting the PDE onto the Laplacian eigenbasis yields decoupled ODEs with analytical solution:
\begin{equation}
    c_k(t) = c_k(0)e^{-D\lambda_k t} + \frac{f_k}{D\lambda_k}\left(1 - e^{-D\lambda_k t}\right),
\end{equation}
where $f_k = \langle f, \phi_k \rangle$ are the forcing coefficients. The network prediction is then constructed as:
\begin{equation}
\label{eq:heat_inhomo_reconstruction}
\hat{u}_i(x,t; \boldsymbol{\theta}) = \sum_{k=1}^{400} \left[ c_{k,i}(0; \boldsymbol{\theta}) e^{-D\lambda_k t} + \frac{f_{k,i}}{D\lambda_k}\left(1 - e^{-D\lambda_k t}\right) \right] \phi_k(x),
\end{equation}
where the branch network learns only the initial coordinates $c_{k,i}(0)$, while the forcing coefficients $f_{k,i} = \langle f_i, \phi_k \rangle$ are pre-computed for each sample based on the corresponding eigenfunctions.

Both the initial condition $u_0(x)$ and the forcing term $f(x)$ were sampled independently from GRFs. Ground truth solutions were computed using the Implicit Euler scheme with $D=0.02$, $T=1.0$, and $dt=0.0025$. Snapshots were saved at 10 uniformly spaced time points.

For FEENet, the branch network takes initial conditions $u_0(x)$ as inputs and learns the initial spectral coordinates $\{c_{k,i}(0)\}$. The forcing coefficients $f_{k,i}$ are pre-computed, while the eigenvalues $\{\lambda_k\}_{k=1}^{400}$ are fixed geometric quantities embedded in the analytical ODE layer. Unlike the previous cases, only $u_0(x)$ is Z-score normalized while $f_k$ and outputs $u(x,t)$ are kept in physical units to preserve exact analytical relations. The dataset is:
\begin{equation}
    \mathcal{T}_{\text{FEENet}} = \left\{ \left( u_{0,i}(x), t_j \right), u_i(x, t_j) \right\}_{i=1,j=1}^{2000,10}.
\end{equation}
The branch network architecture is $[P, 400]$.

For MIONet, standard Z-score normalization is applied to all inputs and outputs. The dataset structure is:
\begin{equation}
    \mathcal{T}_{\text{MIONet}} = \left\{ \left( u_{0,i}(x), f_i(x), (x,t) \right), u_i(x,t) \right\}_{i=1}^{2000}.
\end{equation}
The network employs two branch networks: one for initial conditions $[P, 400, 400, 400]$ and one for forcing terms $[P, 400, 400, 400]$, plus a trunk network $[d+1, 400, 400]$ for spatiotemporal coordinates.

Learning rates are $4\times 10^{-5}$ (Unit Square, Fins) and $2\times 10^{-6}$ (Bunny) for FEENet, and $5\times 10^{-5}$ (Unit Square, Fins) and $1\times 10^{-5}$ (Bunny) for MIONet.

\begin{figure}[htbp]
  \centering
  \includegraphics[width=\linewidth]{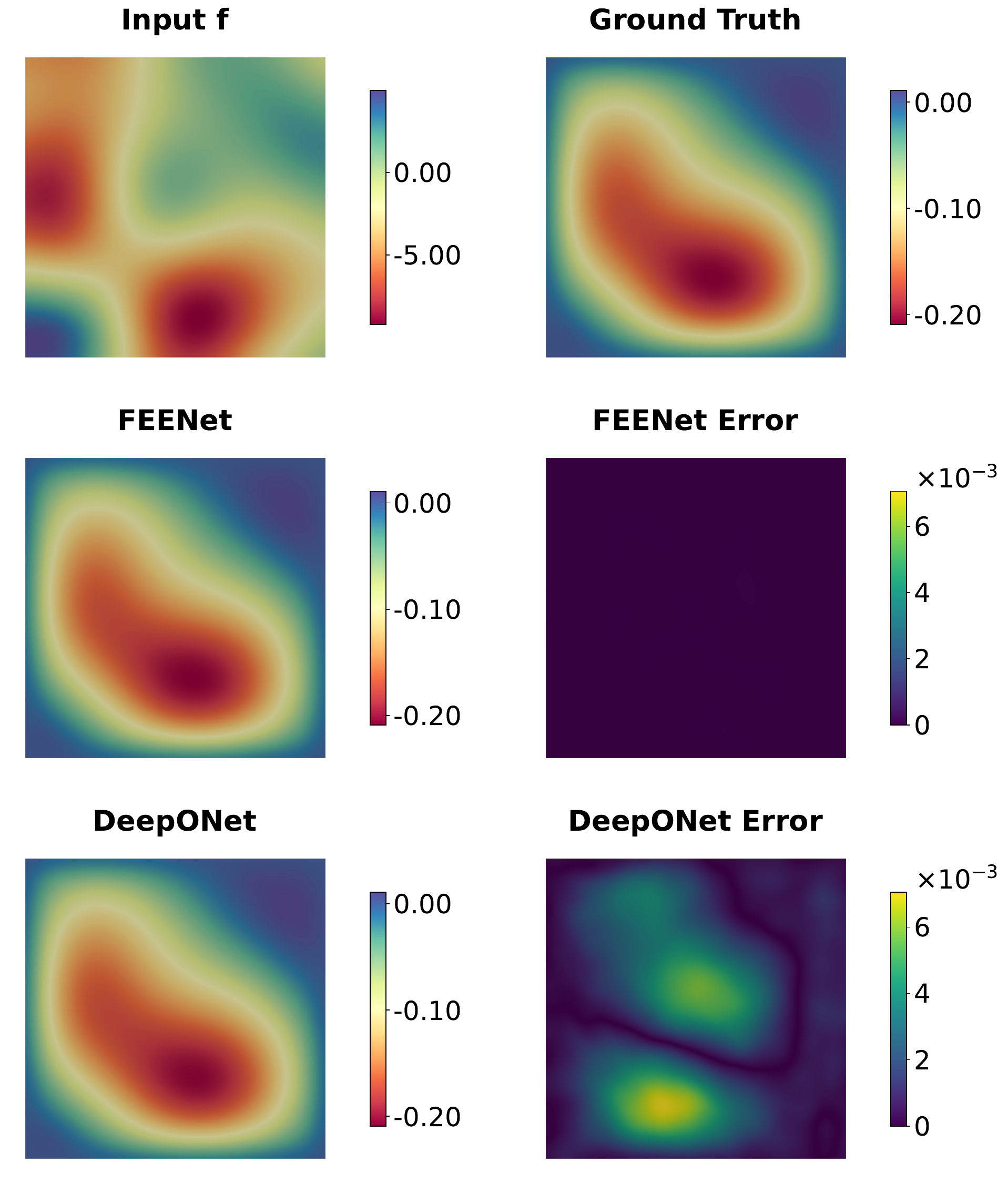}
  \caption{\textbf{Poisson problem on the Unit Square.} Randomly selected representative result for (top row) input function and corresponding reference solution; 
  (middle row) FEENet prediction and the corresponding absolute error. 
  (bottom row) DeepONet prediction and the corresponding absolute error.}
  \label{fig:Poisson_square_model_comparison}
\end{figure}

\begin{figure}[htbp]
  \centering
  \includegraphics[width=\linewidth]{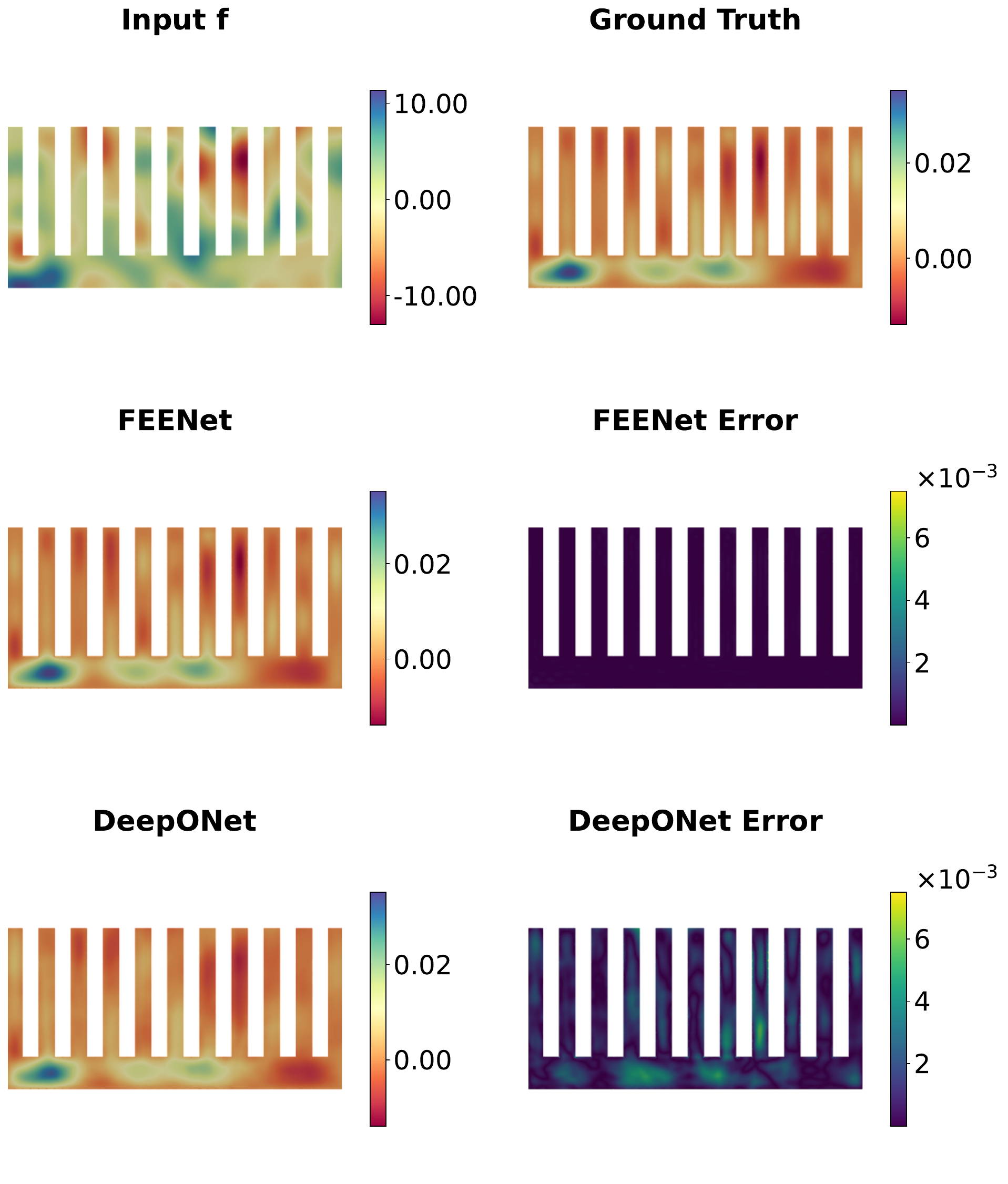}
  \caption{\textbf{Poisson problem on the Fins.} Randomly selected representative result for (top row) input function and corresponding reference solution; 
  (middle row) FEENet prediction and the corresponding absolute error; 
  (bottom row) DeepONet prediction and the corresponding absolute error.}
  \label{fig:Poisson_fins_model_comparison}
\end{figure}

\begin{figure}[htbp]
  \centering
  \includegraphics[width=\linewidth]{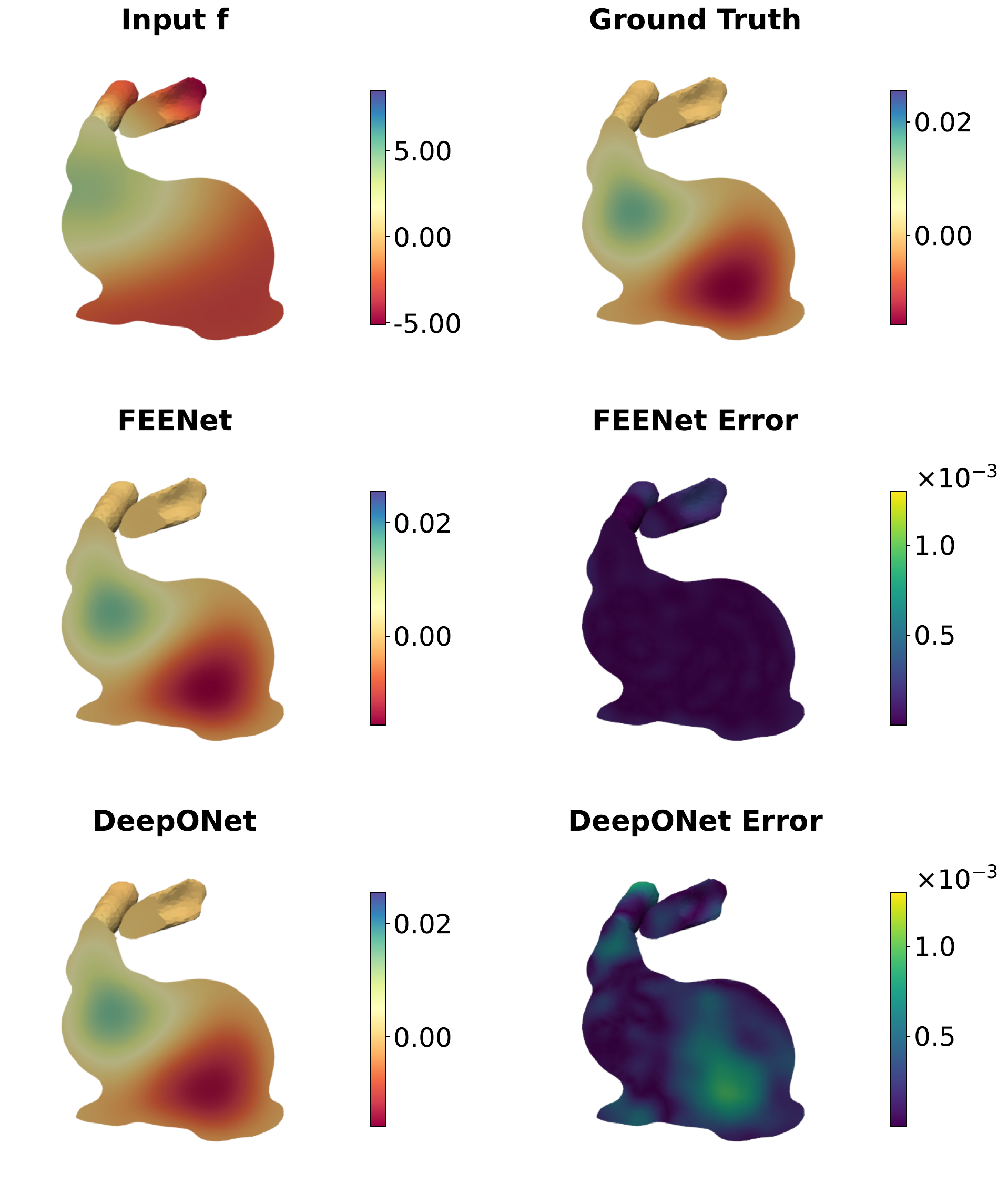}
  \caption{\textbf{Poisson problem on the Bunny.} Randomly selected representative result for
  (top row) input function and corresponding reference solution;
  (middle row) FEENet prediction and the corresponding absolute error;
  (bottom row) DeepONet prediction and the corresponding absolute error.}
  \label{fig:Poisson_bunny_model_comparison}
  \end{figure}

% Hom. heat square
\begin{figure}[htbp]
  \centering
  \includegraphics[width=\linewidth,height=0.90\textheight,keepaspectratio]{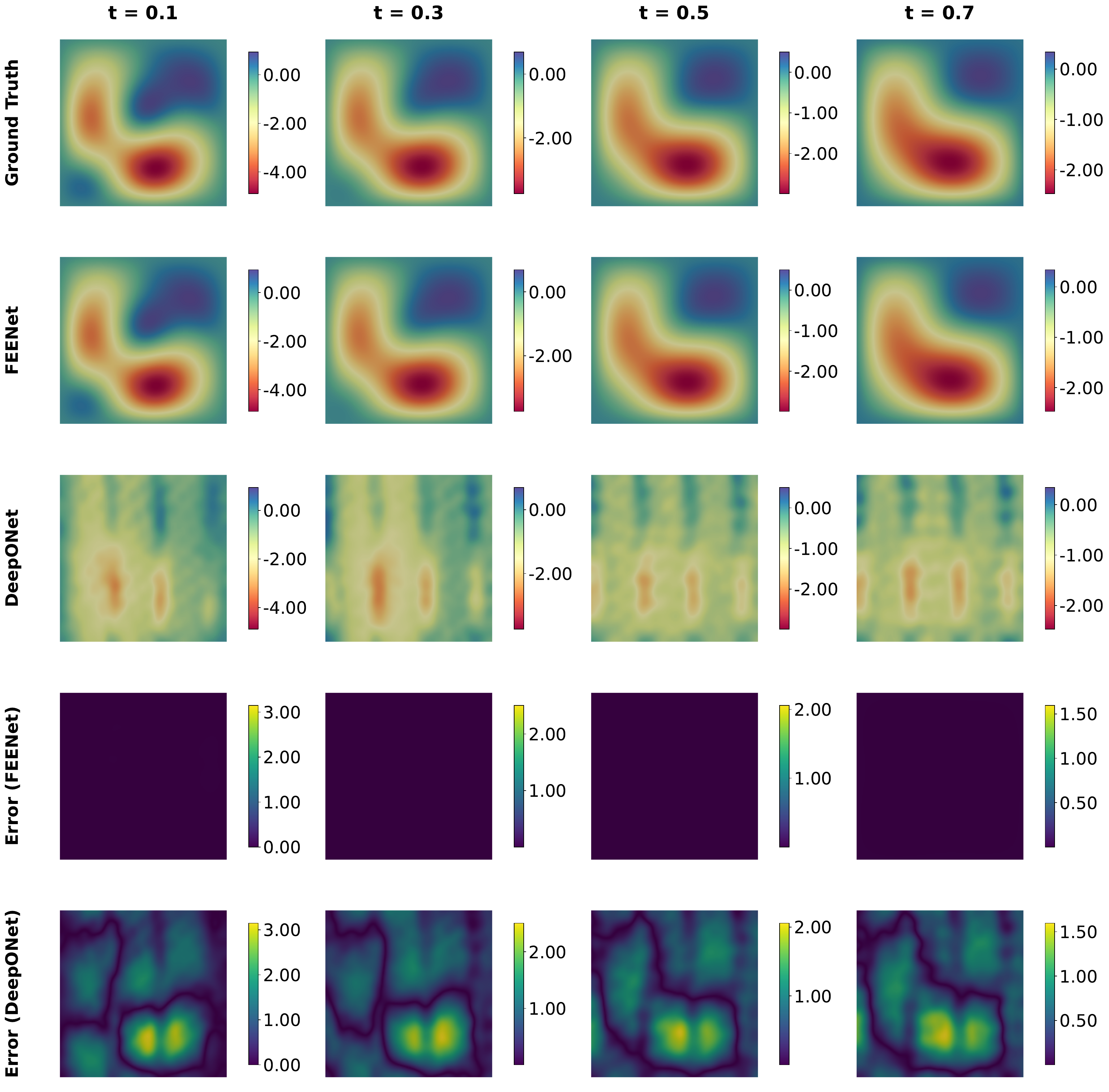}
  \caption{\textbf{Homogeneous heat problem on the Unit Square.}
  From top to bottom: Reference solutions at selected time steps; FEENet predictions; DeepONet predictions; absolute error of FEENet; absolute error of DeepONet.}
  \label{fig:hom_heat_square_model_comparison}
\end{figure}

% Hom. heat fins
\begin{figure}[htbp]
  \centering
  \includegraphics[width=\linewidth,height=0.90\textheight,keepaspectratio]{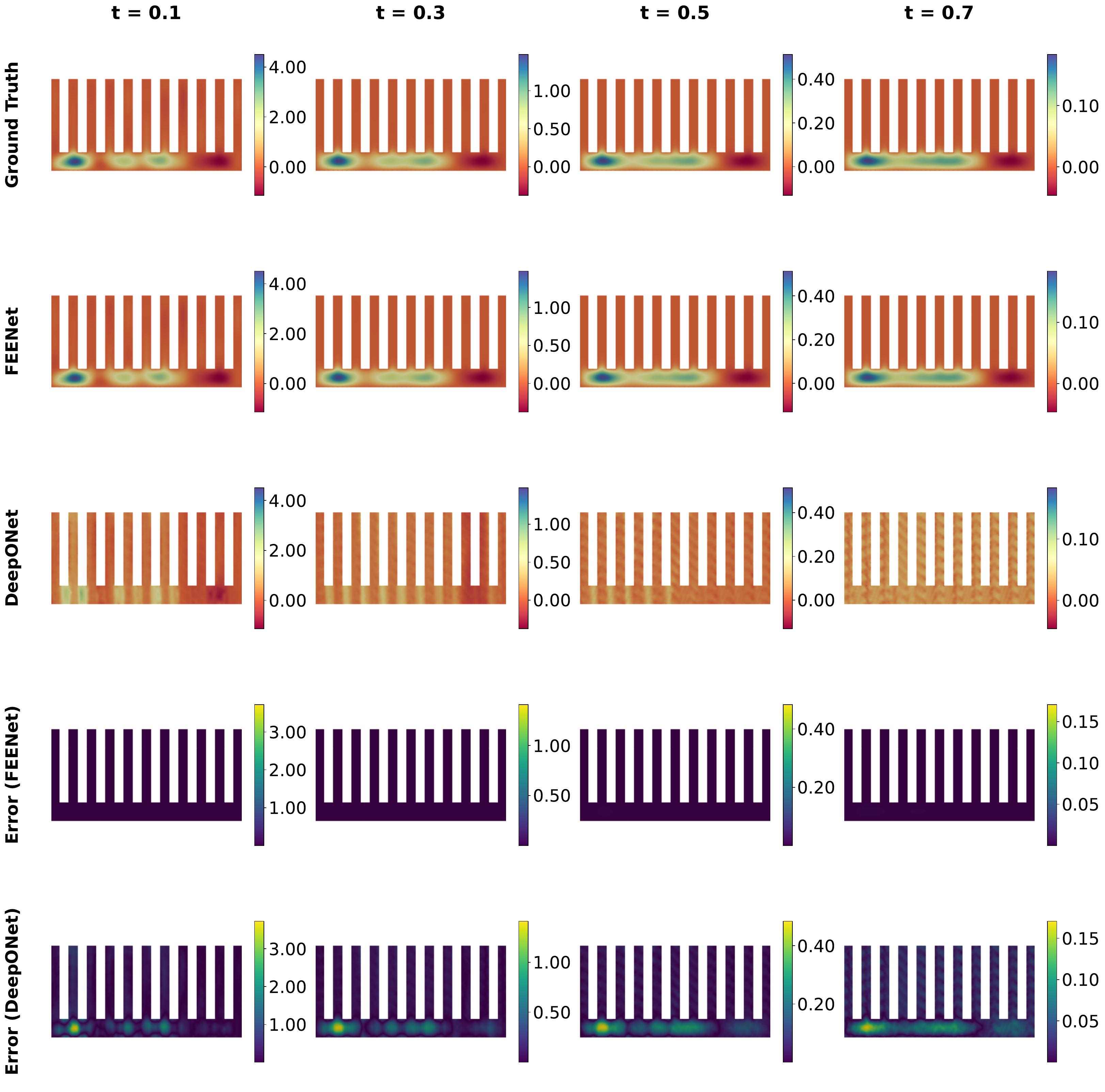}
  \caption{\textbf{Homogeneous heat problem on the Fins.}
  From top to bottom: Reference solutions at selected time steps; FEENet predictions; DeepONet predictions; absolute error of FEENet; absolute error of DeepONet.}  
  \label{fig:hom_heat_fins_model_comparison}
\end{figure}

% Hom. heat bunny
\begin{figure}[htbp]
  \centering
  \includegraphics[width=\linewidth,height=0.90\textheight,keepaspectratio]{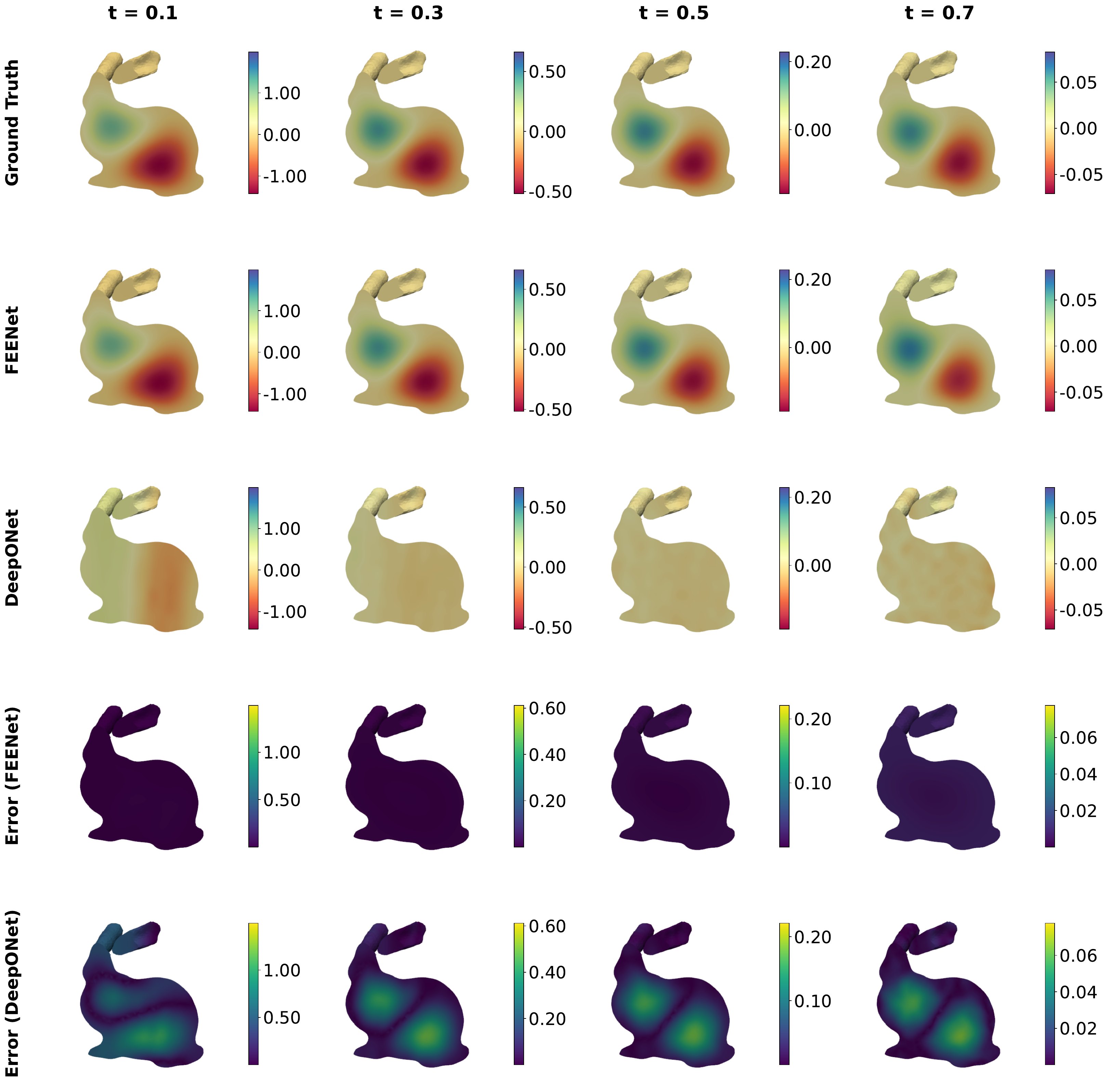}
  \caption{\textbf{Homogeneous heat problem on the Bunny.}
  From top to bottom: Reference solutions at selected time steps; FEENet predictions; DeepONet predictions; absolute error of FEENet; absolute error of DeepONet.}  
  \label{fig:hom_heat_bunny_model_comparison}
\end{figure}

% Inhom. heat square
\begin{figure}[htbp]
  \centering
  \includegraphics[width=\linewidth,height=0.90\textheight,keepaspectratio]{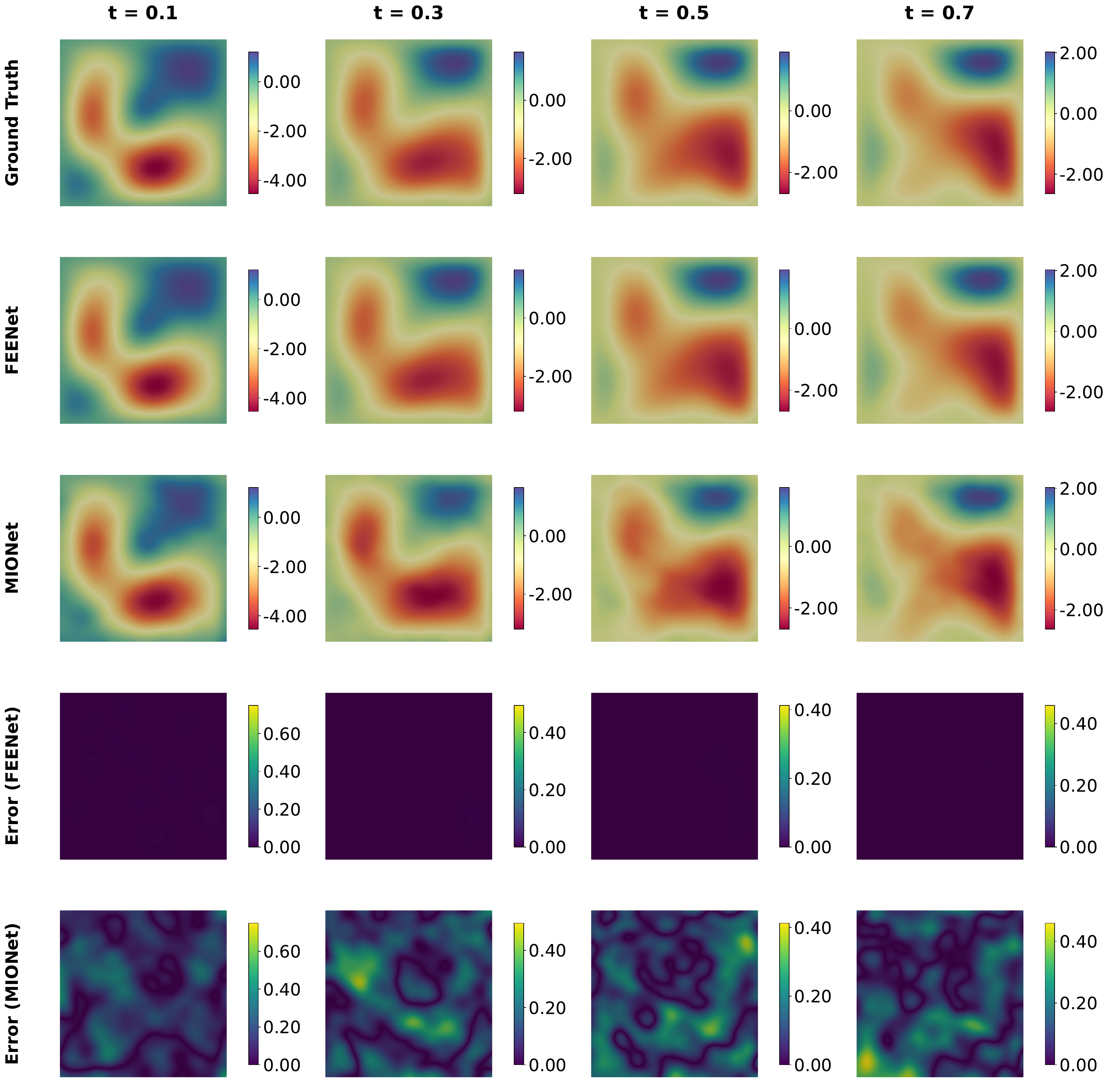}
  \caption{\textbf{Inhomogeneous heat problem on the Unit Square.}
  From top to bottom: Reference solutions at selected time steps; FEENet predictions; DeepONet predictions; absolute error of FEENet; absolute error of DeepONet.}  
  \label{fig:inhom_heat_square_model_comparison}
\end{figure}

% Inhom. heat fins
\begin{figure}[htbp]
  \centering
  \includegraphics[width=\linewidth,height=0.90\textheight,keepaspectratio]{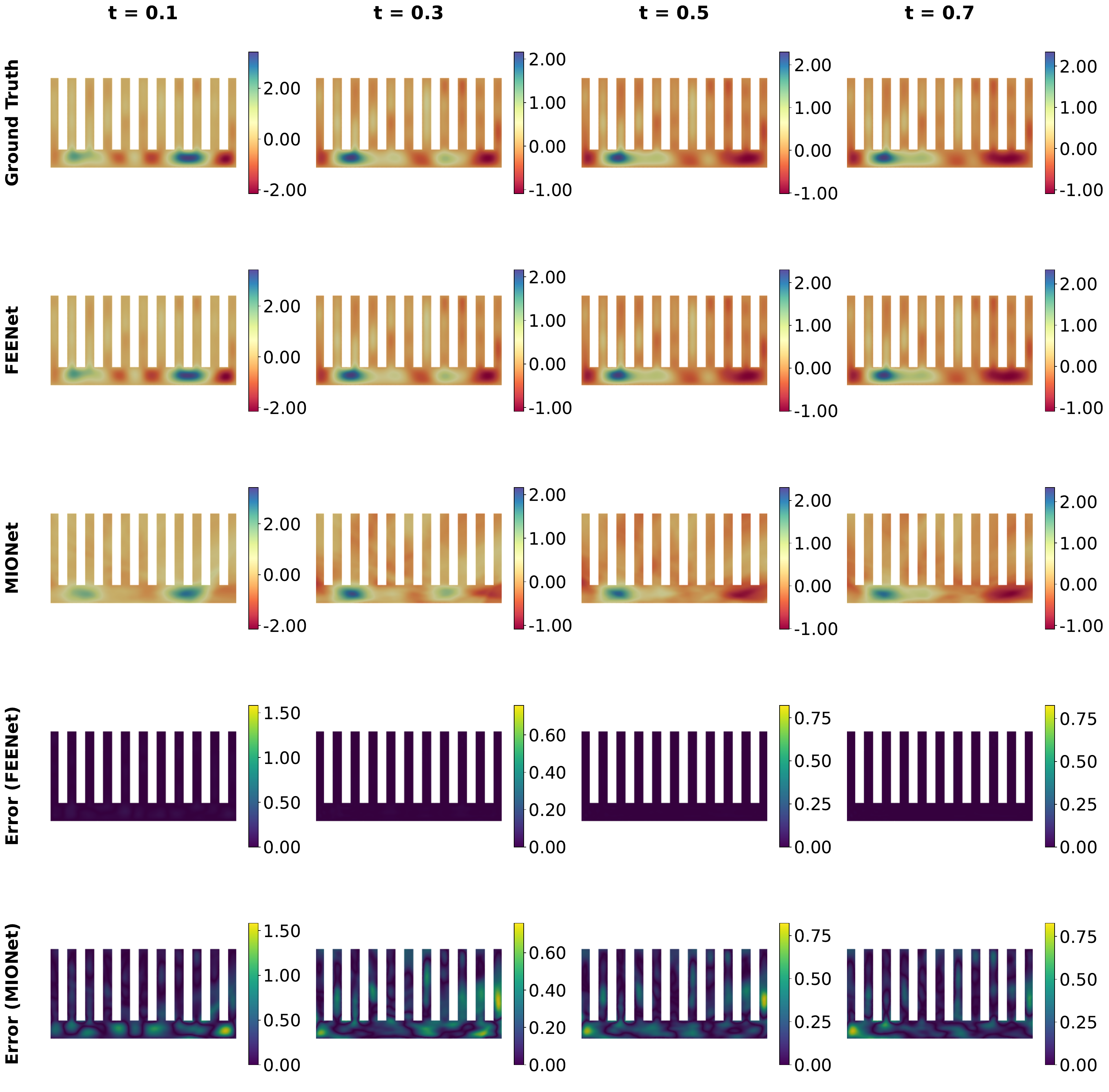}
  \caption{\textbf{Inhomogeneous heat problem on the Fins.}
  From top to bottom: Reference solutions at selected time steps; FEENet predictions; DeepONet predictions; absolute error of FEENet; absolute error of DeepONet.}    
  \label{fig:inhom_heat_fins_model_comparison}
\end{figure}

% Inhom. heat bunny
\begin{figure}[htbp]
  \centering
  \includegraphics[width=\linewidth,height=0.90\textheight,keepaspectratio]{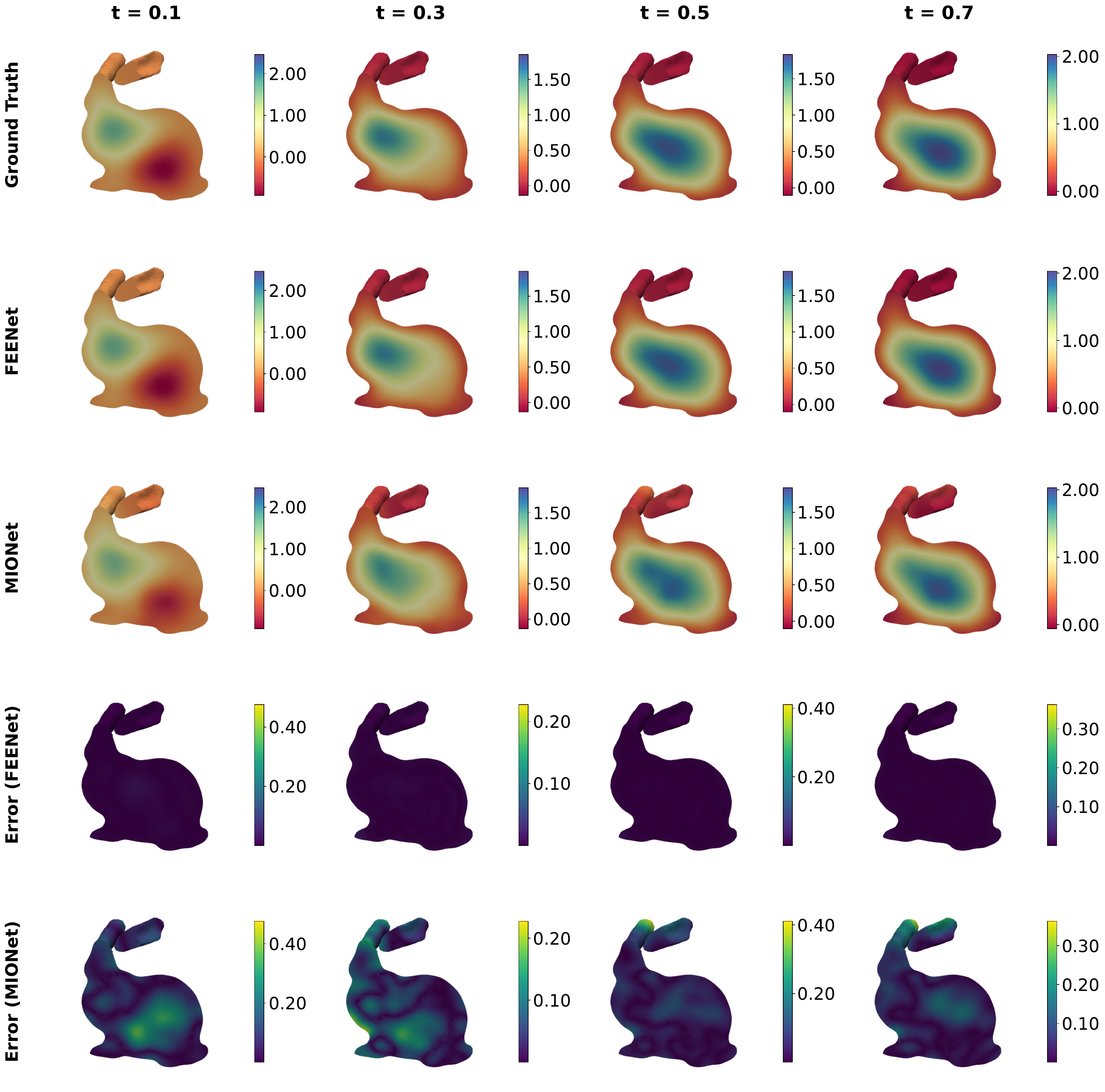}
  \caption{\textbf{Inhomogeneous heat problem on the Bunny.}
  From top to bottom: Reference solutions at selected time steps; FEENet predictions; DeepONet predictions; absolute error of FEENet; absolute error of DeepONet.}    
  \label{fig:inhom_heat_bunny_model_comparison}
\end{figure}

% ========== TABLE: Model training results ==========
% ========== L2 mesh + H1 mesh ==========
\begin{table}[ht]
\centering
\scriptsize
\caption{\textbf{Training cost and accuracy comparison across PDEs and geometries.}
The table reports the number of trainable parameters, wall-clock training time (in minutes), and relative errors in the $L^2$ and $H^1$ norms for different models and geometries. For FEENet, the reported time is given as \emph{training time + eigenbasis generation time}, where the latter is a one-time cost per geometry (see Table~\ref{tab:computation_cost_for_eigenpairs}) and can be amortized across multiple PDE problems on the same domain.}
\label{tab:training}
\setlength{\tabcolsep}{9.5pt} 
\renewcommand{\arraystretch}{1.5}
\begin{tabularx}{\textwidth}{ll l >{\raggedleft\arraybackslash}p{1.5cm} 
>{\raggedleft\arraybackslash}p{1.6cm} 
>{\raggedleft\arraybackslash}p{1.9cm} 
>{\raggedleft\arraybackslash}p{1.9cm}}
\toprule
\textbf{PDE} & \textbf{Domain} & \textbf{Model} & \textbf{Number of param.} & \textbf{Train time (min)} & {\textbf{Relative $L^2$ error}} & \textbf{Relative $H^1$ error} \\
\midrule

% ===== Poisson Equation =====
\multirow{6}{*}{Poisson}
 & \multirow{2}{*}{Square}
   & DeepONet     & 1,133,201  & 12.67 & $1.88\times10^{-2}$ & $4.18\times10^{-2}$ \\
 & & FEENet  & 490,401 & 4.16+0.12  & $5.79\times10^{-4}$ & $4.60\times10^{-3}$ \\
\cmidrule(l){2-7} 

 & \multirow{2}{*}{Fins}
   & DeepONet     & 4,198,801 & 73.65 & $2.13\times10^{-1}$ & $4.86\times10^{-1}$ \\
 & & FEENet  & 3,556,001 & 26.66+1.07 & $1.44\times10^{-2}$ & $6.71\times10^{-2}$ \\
 \cmidrule(l){2-7}

 & \multirow{2}{*}{Bunny}
   & DeepONet     & 16,869,601 & 361.82 & $1.83\times10^{-2}$ & $9.74\times10^{-2}$ \\
 & & FEENet  & 16,226,401 & 146.72+6.52 & $7.98\times10^{-3}$ & $4.26\times10^{-2}$ \\

% ===== Homogeneous Heat Equation =====
\midrule
\multirow{6}{*}{Homog.\ Heat}
 & \multirow{2}{*}{Square}
   & DeepONet     & 1,133,601 & 87.04 & $5.59\times10^{-1}$ & $1.99\times10^{-0}$ \\
 & & FEENet  & 490,400 & 3.71+0.12 & $3.38\times10^{-3}$ & $8.39\times10^{-3}$ \\
 \cmidrule(l){2-7}
 
 & \multirow{2}{*}{Fins}
   & DeepONet     & 4,199,201 & 676.0 & $7.22\times10^{-1}$ & $1.62\times10^{-0}$ \\
 & & FEENet  & 3,556,000 & 24.82+1.07 & $8.29\times10^{-3}$ & $1.07\times10^{-2}$ \\
 \cmidrule(l){2-7}
 
 & \multirow{2}{*}{Bunny}
   & DeepONet     & 16,870,001 & 3092.74 & $6.35\times10^{-1}$ & $1.17\times10^{-0}$ \\
 & & FEENet  & 16,226,400 & 139.58+6.52 & $1.56\times10^{-2}$ & $1.52\times10^{-2}$ \\

% ===== Inhomogeneous Heat Equation =====
\midrule
\multirow{6}{*}{Inhomog.\ Heat}
 & \multirow{2}{*}{Square}
   & MIONet       & 1,784,401 & 62.22 & $1.32\times10^{-1}$ & $2.54\times10^{-1}$ \\
 & & FEENet  & 490,400 & 4.16+0.12 & $1.47\times10^{-3}$ & $7.82\times10^{-3}$ \\
 \cmidrule(l){2-7}
 
 & \multirow{2}{*}{Fins}
   & MIONet       & 7,915,601 & 485.5 & $5.82\times10^{-1}$ & $7.77\times10^{-1}$ \\
 & & FEENet  & 3,556,000 & 26.34+1.07 & $1.96\times10^{-2}$ & $5.95\times10^{-2}$ \\
 \cmidrule(l){2-7}
 
 & \multirow{2}{*}{Bunny}
   & MIONet       & 33,256,801 & 2208.42 & $1.59\times10^{-1}$ & $3.49\times10^{-1}$ \\
 & & FEENet  & 16,226,400 & 143.57+6.52 & $8.57\times10^{-3}$ & $4.01\times10^{-2}$ \\

\bottomrule
\end{tabularx}
\end{table}

\section{Discussion}
\label{sec:discussion}
As a hybrid methodology, FEENet is designed to inherit complementary strengths from both classical numerical analysis and modern operator learning. On the one hand, by constructing its representation space from finite element eigenfunctions, FEENet retains the geometric fidelity, stability, and approximation guarantees characteristic of FEM-based discretizations on complex domains. On the other hand, by learning the solution operator in this fixed, operator-adapted basis, FEENet achieves the computational efficiency, rapid evaluation, and many-query capability typical of neural operator models. We aim to highlight and quantify some of these features in this section.

\subsection{Resolution Independence Study}
A defining characteristic of neural operators is resolution independence—the capability to evaluate the learned solution at arbitrary spatial coordinates without model retraining. Unlike traditional finite element methods that rely on fixed mesh connectivity, FEENet learns the mapping from the input function space to the spectral coordinate space. Since the solution is represented as a linear combination of continuous eigenfunctions:
\begin{equation}
    u(x) \approx \sum_{k=1}^M c_k \phi_k(x),
\end{equation}
the spatial resolution of the output is decoupled from the training discretization and is determined solely by the evaluation of the basis functions $\phi_k(x)$. Once the eigenpairs are computed, the solution can be queried at any location $x \in \Omega$.

To empirically validate this property, we conducted a resolution independence study using the Poisson equation on the Fins geometry. The model was trained using a dataset generated with a low density (totaling 8,889 points). During the inference phase, with fixed model weights, we evaluated prediction accuracy on three kinds of inquiry structures with increasing discretization densities. In all cases, the inquiry points are generated from uniform, equally spaced structured grids covering the domain.

\begin{enumerate}
    \item \textbf{Low Density (Training Resolution):} 
    The structured sampling grid used during training (8,889 sensor points), 
    serving as the only resolution seen by the model during optimization.

    \item \textbf{Medium Density:} 
    A denser structured inquiry grid that was \emph{not} used during training, 
    containing approximately four times as many sensor points as the training 
    resolution (33,938 points).

    \item \textbf{High Density:} 
    A significantly finer inquiry grid with over fifteen times the number of 
    sensor points compared to training (135,175 points), used to assess 
    prediction consistency at high spatial resolution.
\end{enumerate}

\begin{figure}[htbp]
    \centering
    \includegraphics[width=\linewidth]{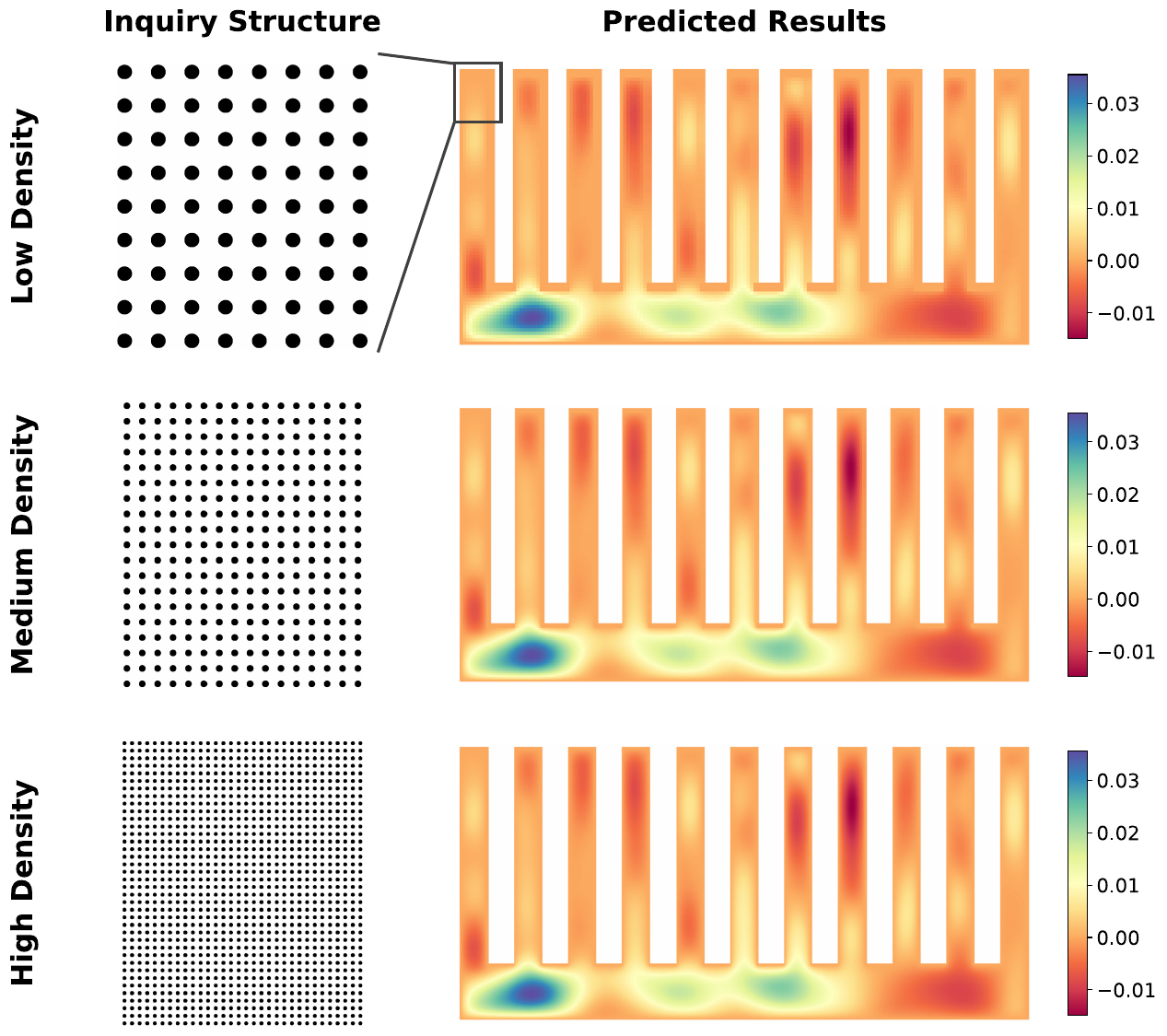}
    \caption{\textbf{Resolution independence in FEENet.} We use the Poisson problem on the Fins geometry
    with three sampling densities from top to bottom. The left column shows the \textit{inquiry grid density} as a zoomed-in view of the structured sampling grids at each density level.
    Right column depicts the corresponding FEENet prediction.}
    \label{fig:resolution_vis}
\end{figure}

As shown in Table \ref{tab:resolution_independence}, the relative $L^2$ and $H^1$ errors remain remarkably consistent across all resolutions. Despite the number of evaluation points increasing by more than an order of magnitude, the error metrics show negligible deviation. This stability confirms that FEENet effectively learns the continuous solution operator and avoids overfitting to the specific spatial discretization of the training data.

\begin{table}[ht]
\centering
\small
\caption{\textbf{Quantitative evaluation of resolution independence.}
The FEENet model was trained exclusively on the \textit{low density} (training) resolution.
The relative $L^2$ and $H^1$ errors were subsequently evaluated on unseen inquiry grids of \textit{medium} and \textit{high density} without retraining.
The closeness of errors across all densities demonstrate evidence that the learned solution is invariant to the grid resolution.}
\label{tab:resolution_independence}
\begin{tabular}{l c c c}
\toprule
\textbf{Discretization density} & \textbf{Number of points} & \textbf{Relative $L^2$ error} & \textbf{Relative $H^1$ error} \\
\midrule
Low Density (Training)  & 8,889   & $1.45 \times 10^{-2}$ & $6.71 \times 10^{-2}$ \\
Medium Density          & 33,938  & $1.59 \times 10^{-2}$ & $7.08 \times 10^{-2}$ \\
High Density            & 135,175 & $1.60 \times 10^{-2}$ & $7.48 \times 10^{-2}$ \\
\bottomrule
\end{tabular}
\end{table}

It is worth noting that while the FEENet inference is resolution-independent—allowing queries at arbitrary coordinates—the ultimate accuracy of the predicted solution depends on the quality of the spectral basis. Since the eigenfunctions $\phi_k(x)$ are pre-computed using a finite element solver, they inherently carry the discretization error associated with the offline FEM mesh. If the initial FEM mesh is too coarse the basis functions fail to resolve high-frequency geometric features and the error will subsequently propagate to the FEENet predictions.

\subsection{Effect of the Number of Eigenfunctions}
\label{eigen:study}

We investigate the influence of the truncation number of eigenfunctions $M$ on the accuracy of the proposed FEENet. We again recourse to the Poisson problem on the fin geometry for this numerical experiment. To this end, we pre-computed the first $800$ Laplacian eigenfunctions of the Fins geometry. We then trained eight separate models by varying $M$ from 100 to 800 to assess performance. In this process, we use the identical PDE data as described in Section~\ref{poisson_problem} but the dimension of the outputs of branch net in FEENet varies to match the numbers of eigenfunctions being used. The spectral coordinates were learned using a single fully-connected layer mapping the input sensors to the $M$ coefficients (layer size $[8889, M]$). All models were consistently trained with a batch size of $256$, a learning rate of $4\times10^{-5}$, and $100{,}000$ iterations.

\begin{figure}[htbp]
    \centering
    \includegraphics[width=0.85\textwidth]{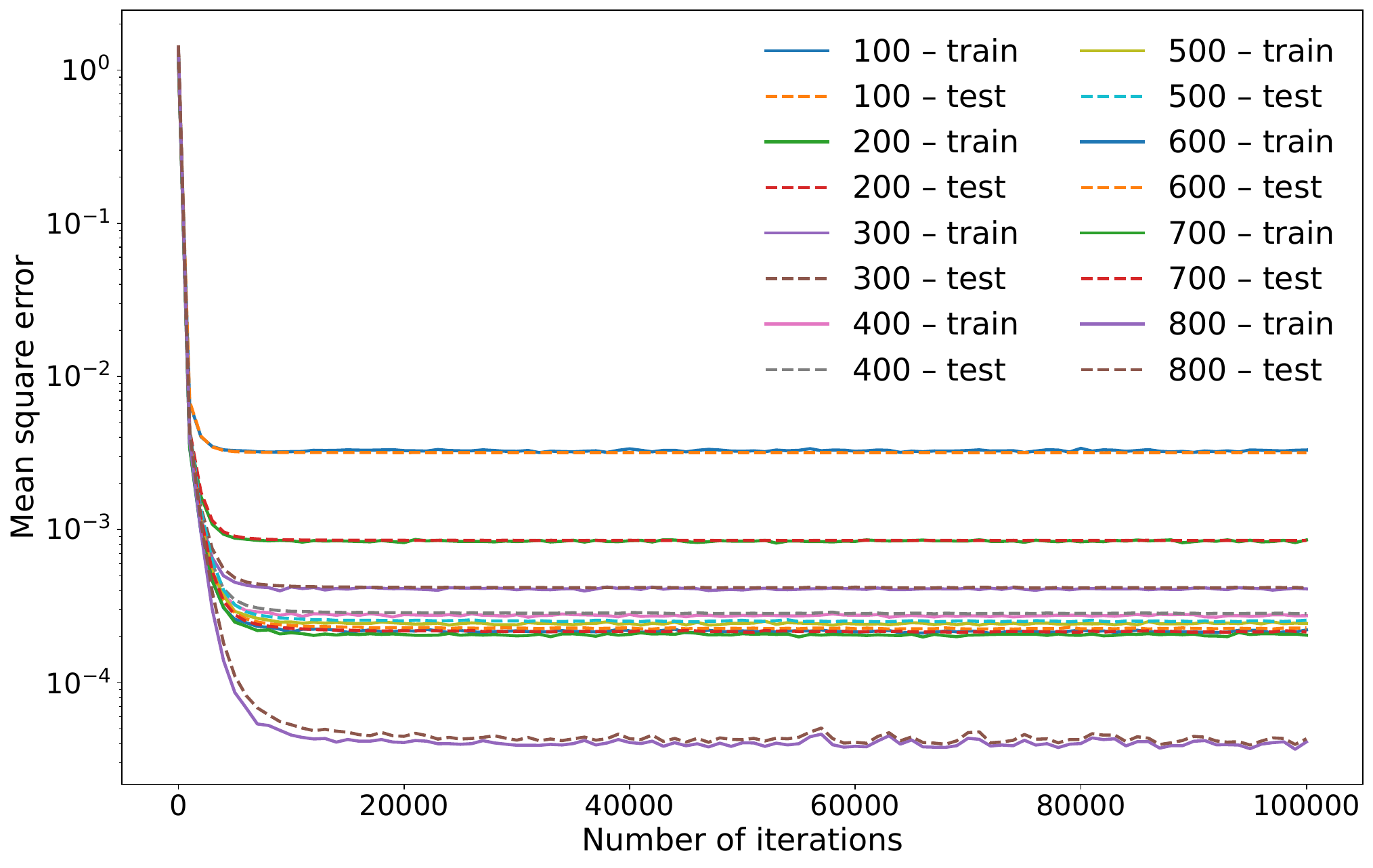}
    \caption{\textbf{Effect of the number of eigenfunctions on training and testing errors for the Poisson problem on the Fins geometry.}
    The curves show the mean squared error for solving the Poisson equation as the number of eigenfunctions $M$ increases.
    The errors decrease monotonically as $M$ increases from 100 to 800, with a particularly marked reduction observed at $M = 800$, indicating the contribution of higher-frequency modes.}
    \label{fig:metric}
\end{figure}

\begin{table}[htbp]
\renewcommand{\arraystretch}{1.4}
\centering

\caption{\textbf{Performance comparison with varying spectral resolution $M$.} The table details the model complexity (parameters), computational cost (training time), and prediction accuracy for the Poisson problem. While training time grows linearly with $M$, the error decreases, particularly at $M=800$.}

\resizebox{\textwidth}{!}{%
\begin{tabular}{c c c c c} % Adjusted column count matching data
\toprule
\textbf{\# Eigenfunctions (M)} & \textbf{\# Parameters} & \textbf{Train Time (min)} & \textbf{Relative $L^2$ error} & \textbf{Relative $H^1$ error} \\
\midrule
100 & 889,001  & 11.58 & $5.33\times10^{-2}$ & $1.41\times10^{-1}$ \\
200 & 1,778,001 & 17.20 & $2.65\times10^{-2}$ & $9.63\times10^{-2}$ \\
300 & 2,667,001 & 22.87 & $1.80\times10^{-2}$ & $7.57\times10^{-2}$ \\
400 & 3,556,001 & 27.58 & $1.44\times10^{-2}$ & $6.71\times10^{-2}$ \\
500 & 4,445,001 & 35.91 & $1.35\times10^{-2}$ & $6.41\times10^{-2}$ \\
600 & 5,334,001 & 39.91 & $1.28\times10^{-2}$ & $6.13\times10^{-2}$ \\
700 & 6,223,001 & 44.07 & $1.25\times10^{-2}$ & $5.99\times10^{-2}$ \\
800 & 7,112,001 & 50.83 & $\mathbf{4.89\times10^{-3}}$ & $\mathbf{2.89\times10^{-3}}$ \\
\bottomrule
\end{tabular}
}
\end{table}

\subsection{FEENet and nonlocal operators}
An additional advantage of FEENet is its natural compatibility with nonlocal operators defined as functions of elliptic operators. Let $\mathcal{L}$ be a self-adjoint strongly elliptic differential operator on a bounded domain $\Omega$ with Dirichlet boundary conditions. By the spectral theorem, for any sufficiently regular function $g:\mathbb{R}^+ \to \mathbb{R}$, one may define the operator $g(\mathcal{L})$ via functional calculus, including important examples such as fractional powers $\mathcal{L}^\alpha$, trigonometric functions $\sin(\mathcal{L})$, and other nonlocal operators. In particular, if $\{\lambda_j, \psi_j\}_{j\ge 1}$ denote the eigenpairs of $\mathcal{L}$, then for any $u\in L^2(\Omega)$ one has
\[
g(\mathcal{L})u = \sum_{j=1}^\infty g(\lambda_j)\,\langle u,\psi_j\rangle\,\psi_j,
\]
with convergence in $L^2(\Omega)$ and, under additional regularity assumptions, in appropriate Sobolev norms.

Since FEENet represents solutions in the eigenfunction basis of $\mathcal{L}$, the action of $g(\mathcal{L})$ is diagonalized, reducing the application or learning of such operators to the evaluation (or regression) of the scalar function $g$ on the spectrum $\{\lambda_j\}$. This enables efficient and interpretable learning of nonlocal operators and operator-valued constitutive laws, avoiding black-box approximation of nonlocal effects in physical space.

\section{Conclusion}
\label{sec:conclusion}
In this work, we propose the Finite Element Eigenfunction Network (FEENet), a novel hybrid framework that synergizes the geometric flexibility of classical Finite Element Methods (FEM) with the approximation power of neural operators. By explicitly incorporating the eigenfunctions of the differential operators as a domain-intrinsic basis, FEENet effectively decouples the challenge of geometric complexity from the operator learning process.

Our numerical experiments, spanning the Poisson problem, the homogeneous heat problem, and the inhomogeneous heat problem with a time-independent forcing term across diverse geometries in 2D and 3D demonstrate the superiority of FEENet over the standard DeepONet and its multi-input variant MIONet. The key findings are summarized as follows:
\begin{enumerate}
    \item \textbf{Geometric Generalization:} FEENet seamlessly handles irregular non-convex domains without requiring complex coordinate transformations or domain decomposition.
    \item \textbf{Efficiency and Convergence:} By offloading the geometric encoding to a pre-computed spectral basis, the learning task is reduced to determining spectral coordinates. This results in very fast convergence  and reduced training costs.
    \item \textbf{Physical Consistency:} For time-dependent problems, FEENet integrates the analytical temporal evolution of eigenmodes directly into the architecture. This ensures strict adherence to the underlying physics and enhances extrapolation capabilities.
    \item \textbf{Spectral Scalability:} As shown in Section~\ref{eigen:study}, increasing the spectral resolution systematically improves prediction accuracy. 
\end{enumerate}

Furthermore, the spectral architecture of FEENet naturally extends to nonlocal operators defined as functions of eigenfunctions,  where a \textit{nested learning} paradigm can potentially enable approximating such nonlocal operators and pseudo-differential operators \cite{taylor1996partial,di2012hitchhiker}.
FEENet enables this capability on arbitrary geometries  without modification to the network architecture, a task that remains computationally prohibitive for traditional mesh-based fractional solvers.

While this work demonstrates the efficacy of FEENet for problems with homogeneous Dirichlet boundary conditions, extending the framework to more general boundary settings is envisioned and remains an important direction. We envision that FEENet, building on the inspiring developments in neural operator theory such as \cite{lu2019deeponet}, will provide a promising step towards interpretable, geometry-aware scientific machine learning for real-world engineering applications.

\section*{Code Availability}
The source code for the proposed FEENet framework is available at \url{https://github.com/l-shiyuan/FEENet}.

\section*{CRediT authorship contribution statement}
\textbf{Shiyuan Li:} Conceptualization, Methodology, Software, Validation, Formal analysis, Investigation, Data Curation, Writing - Original Draft, Visualization. 
\textbf{Hossein Salahshoor:} Conceptualization, Methodology, Resources, Writing - Review \& Editing, Supervision, Project administration, Funding acquisition.

\section*{Data Availability}
Part of the data are available at \url{https://github.com/l-shiyuan/FEENet}. Additional datasets are available from the corresponding author upon reasonable request.

\section*{Declaration of competing interest}
The authors declare that they have no known competing financial interests or personal relationships that could have appeared to influence the work reported in this paper.

\section*{Acknowledgments}
The first author expresses his gratitude to Professor Thomas P. Witelski of the Department of Mathematics, Duke University, for helpful discussions. This work was partly supported by internal grants from Duke University.

\bibliographystyle{elsarticle-num} 
\bibliography{refs}

\clearpage
\appendix
\section{Laplacian Eigenfunction Visualizations}

This appendix provides visual representations of the Laplacian eigenfunctions computed on the three geometries considered in this work: the Unit Square (Figure~\ref{fig:app_square_lap}), the Fins (Figure~\ref{fig:app_fins_lap}), and the Bunny (Figure~\ref{fig:app_bunny_lap}). These eigenfunctions serve as the fixed spectral basis in the FEENet architecture and encode the geometric structure of each domain.

\subsection{Computational Setting}

For each geometry $\Omega$, the eigenfunctions $\{\phi_k\}_{k=1}^M$ are obtained by solving the Laplacian eigenvalue problem with homogeneous Dirichlet boundary conditions:
\begin{equation}
\begin{cases}
-\Delta \phi_k = \lambda_k \phi_k & \text{in } \Omega, \\
\phi_k = 0 & \text{on } \partial\Omega,
\end{cases}
\label{eq:app_eigen_problem}
\end{equation}
where $\lambda_k$ denotes the $k$-th eigenvalue and $\phi_k$ the corresponding eigenfunction.
The eigenpairs are computed using the finite element method via the DOLFINx library~\cite{baratta2023dolfinx}.

\subsection{Geometric and Spectral Properties}

\begin{figure}[p]
  \centering
  \includegraphics[width=\linewidth,height=0.90\textheight,keepaspectratio]{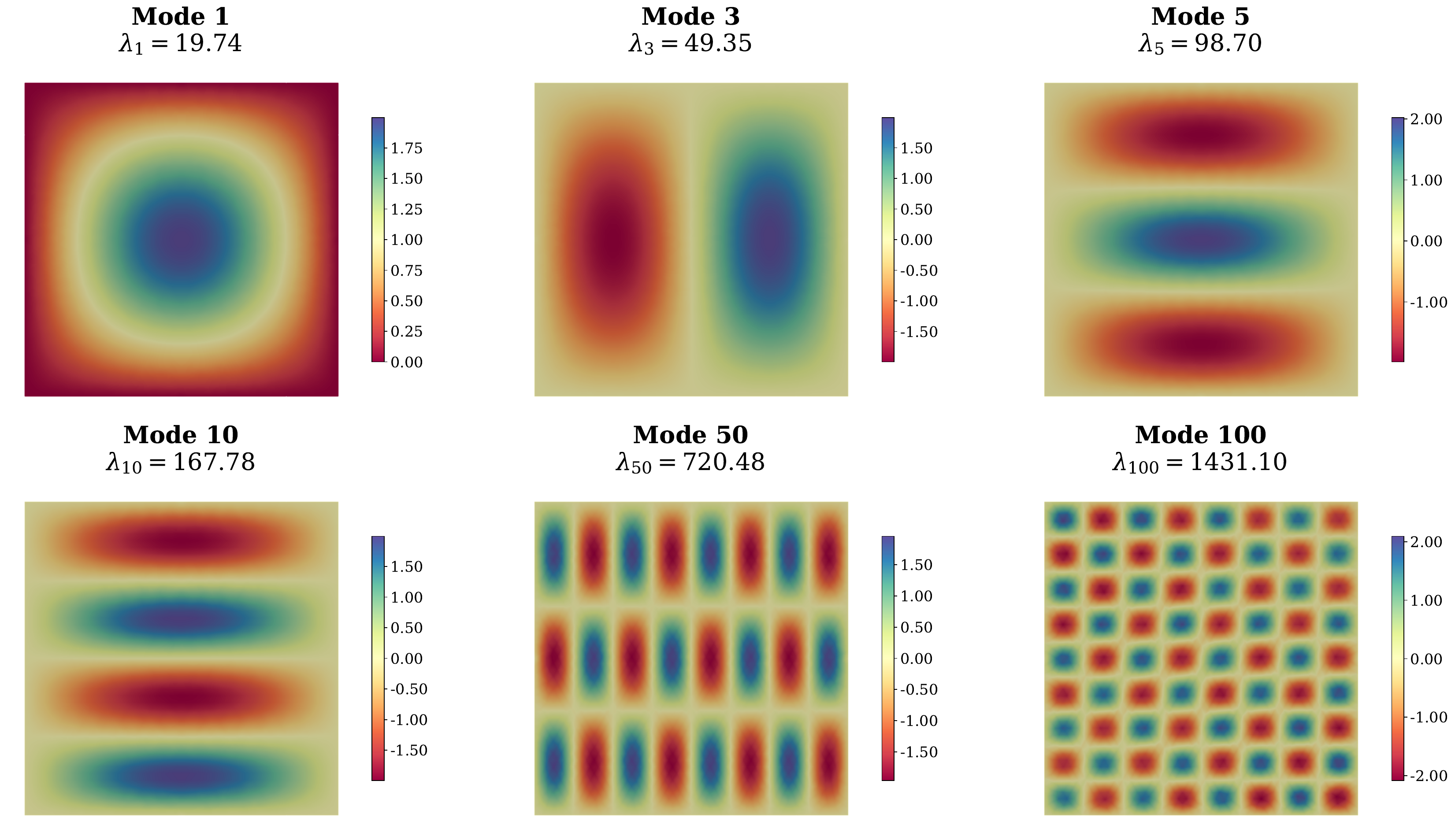}
  \caption{Representative Laplacian eigenfunctions on the Unit Square, with eigenvalues indicated.}
  \label{fig:app_square_lap}
\end{figure}

\begin{figure}[p]
  \centering
  \includegraphics[width=\linewidth,height=0.90\textheight,keepaspectratio]{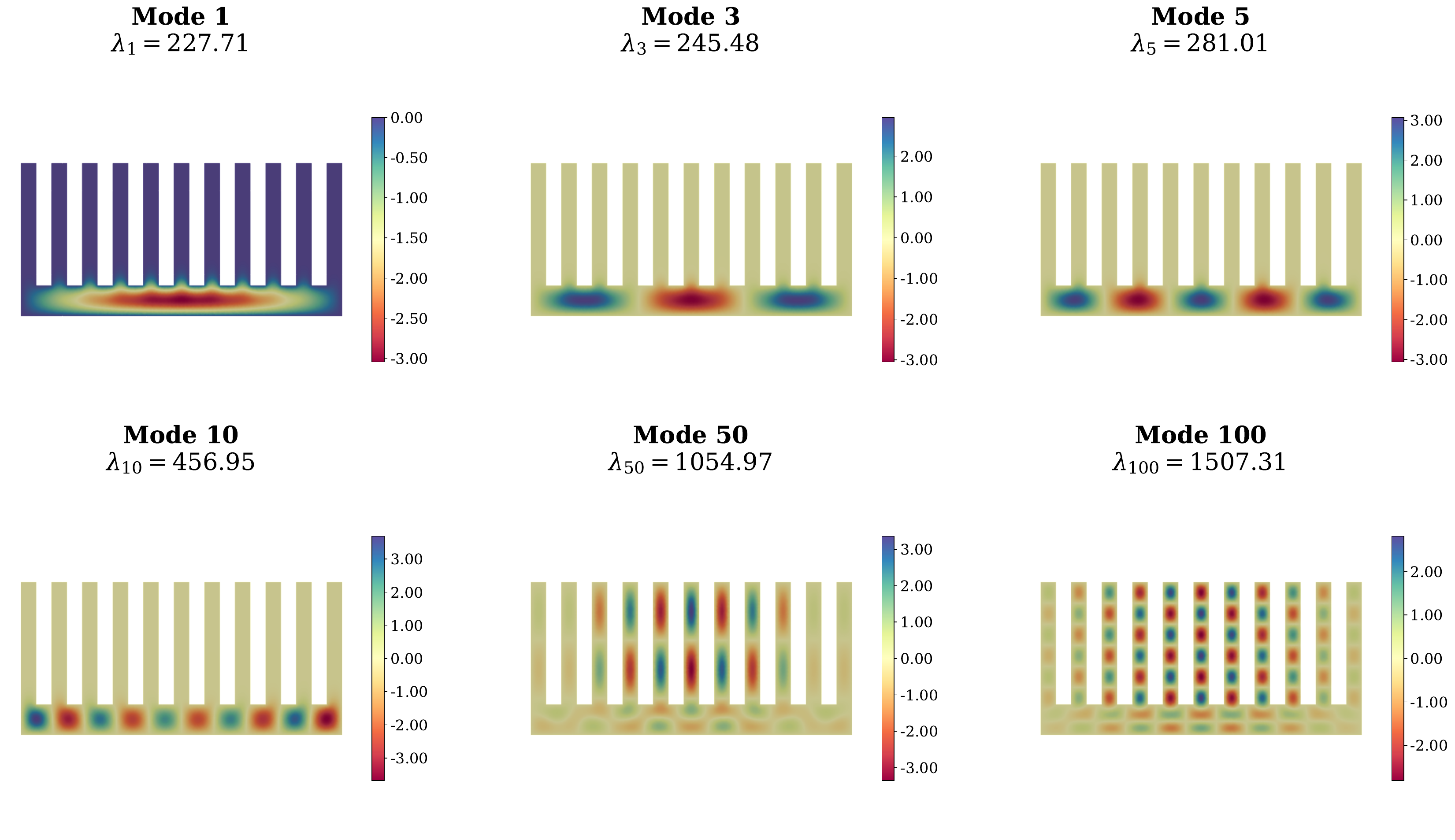}
  \caption{Representative Laplacian eigenfunctions on the Fins, with eigenvalues indicated.}
  \label{fig:app_fins_lap}
\end{figure}

\begin{figure}[p]
  \centering
  \includegraphics[width=\linewidth,height=0.90\textheight,keepaspectratio]{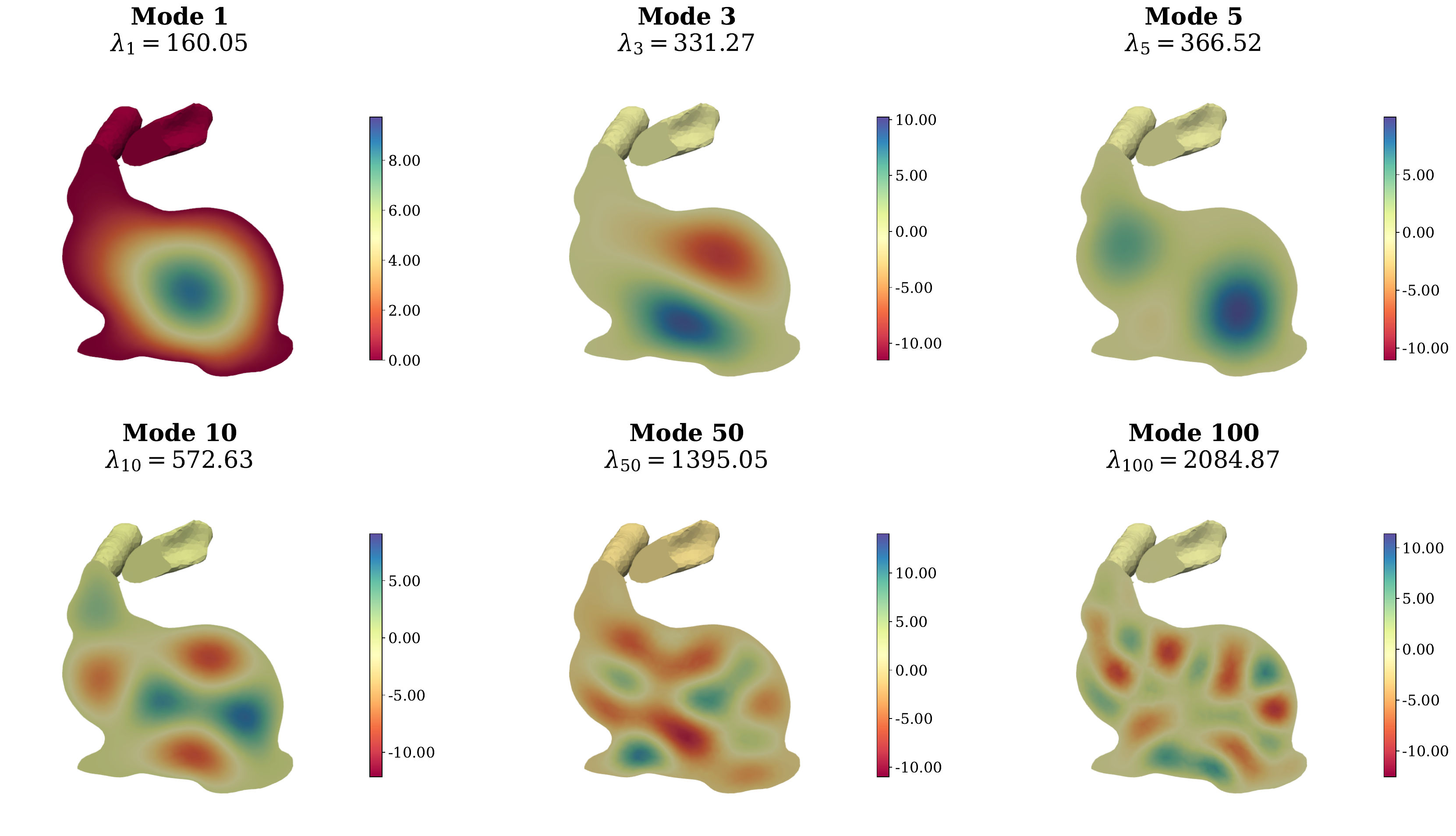}
  \caption{Representative Laplacian eigenfunctions on the Bunny, with eigenvalues indicated.}
  \label{fig:app_bunny_lap}
\end{figure}

The visualizations presented in Figures~\ref{fig:app_square_lap}--\ref{fig:app_bunny_lap} illustrate representative eigenfunctions ordered by increasing eigenvalue. Each eigenfunction exhibits a distinct spatial oscillatory structure, with higher-index eigenfunctions typically displaying more complex nodal patterns. Low-frequency modes with small \(\lambda_k\) capture global geometric features, while high-frequency modes with large \(\lambda_k\) resolve fine-scale details.

For the Unit Square (Figure~\ref{fig:app_square_lap}), the eigenfunctions are products of one-dimensional sine functions and admit the analytical form
\[
\phi_{m,n}(x,y) = \sin(m\pi x)\sin(n\pi y),
\]
with corresponding eigenvalues
\[
\lambda_{m,n} = \pi^2 (m^2 + n^2),
\]
for positive integers $m,n$.
The FEM-computed eigenfunctions closely match this analytical solution, providing a verification of the numerical procedure.

For the Fins geometry (Figure~\ref{fig:app_fins_lap}), the boundary irregularities break the perfect symmetry of the Unit Square geometry.
While the low-frequency eigenfunctions still retain some symmetry, higher-frequency modes exhibit pronounced asymmetry.
In addition, localized oscillations near the thin rectangular protrusions become increasingly prominent in the high-frequency regime, reflecting the strong influence of fine geometric features.

For the three-dimensional Bunny (Figure~\ref{fig:app_bunny_lap}), the eigenfunctions reflect the underlying geometric structure, with increasing spatial complexity at higher frequencies.

% Square — Laplacian

% fins — Laplacian

% Bunny — Laplacian

\section{Gaussian Random Field Generation}
\label{app:Gaussian_generation}

This appendix provides detailed information on the generation of training and testing datasets for the benchmark problems presented in Section~\ref{sec:experiments}. All input functions—including forcing terms for the Poisson and inhomogeneous heat problems, and initial conditions for the heat problems—are sampled from Gaussian Random Fields (GRFs) to ensure diverse and physically meaningful test cases.

\subsection{Covariance Model}

All Gaussian Random Fields are defined by a zero-mean Gaussian distribution
with a stationary Gaussian (squared-exponential) covariance kernel.
Given two spatial locations $\mathbf{x}_1, \mathbf{x}_2 \in \Omega$,
the covariance function is defined as
\begin{equation}
k(\mathbf{x}_1, \mathbf{x}_2)
= \sigma^2
\exp\left(
- \frac{\pi \|\mathbf{x}_1 - \mathbf{x}_2\|^2}{4 \ell^2}
\right),
\end{equation}
where $\sigma^2$ denotes the variance
and $\ell$ is the correlation length scale.
Throughout all experiments, the variance is fixed to $\sigma^2 = 15$.
The correlation length scale $\ell$ is chosen depending on the geometry
to ensure sufficient spatial variability relative to the domain size,
as specified in Section~\ref{sec:experiments}.

\subsection{Randomization Method}

To generate GRF realizations, we employ the Randomization method \citep{hesse2014generating}, a spectral synthesis technique implemented in GSTools \citep{muller2022gstools}. This method approximates the random field as a superposition of Fourier modes with random frequencies and phases:\begin{equation}u(\mathbf{x}) \approx \sqrt{\frac{\sigma^2}{M}} \sum_{m=1}^{M} \left[ Z_{1,m} \cos(\mathbf{k}_m \cdot \mathbf{x}) + Z_{2,m} \sin(\mathbf{k}_m \cdot \mathbf{x}) \right],\label{eq:app_randmeth}\end{equation}where $M$ is the number of modes, $Z_{1,m}, Z_{2,m} \sim \mathcal{N}(0, 1)$ are independent standard normal coefficients, and $\mathbf{k}_m$ are random wave vectors drawn from the normalized spectrum
\begin{equation}
E(\mathbf{k}) = \frac{S(\mathbf{k})}{\sigma^2},
\end{equation}
which defines a probability density function. Here $S(\mathbf{k})$ denotes the power spectral density, i.e., the Fourier transform of the stationary covariance kernel. The computational complexity of this method is $\mathcal{O}(MN)$ for $N$ evaluation points, which scales linearly in $N$. This makes it highly efficient for high-resolution 3D geometries compared to decomposition-based methods (such as Cholesky decomposition) which typically scale as $\mathcal{O}(N^3)$.

\subsection{Boundary Condition Enforcement}

For time-dependent problems requiring homogeneous Dirichlet boundary conditions on initial conditions, raw GRF realizations $u_{\text{GRF}}(x)$ are post-processed using a smooth cutoff function $d(x)$ defined as the solution to:
\begin{equation}
-\Delta d(x) = 1, \quad x \in \Omega,
\qquad
d(x) = 0, \quad x \in \partial \Omega,
\end{equation}
which is subsequently normalized such that
\begin{equation}
\max_{x \in \Omega} d(x) = 1.
\end{equation}
The final admissible initial condition is then constructed as
\begin{equation}
u_{\text{valid}}(x) = u_{\text{GRF}}(x)\, d(x).
\end{equation}
This procedure ensures that $u_{\text{valid}}(x)$ satisfies the prescribed
homogeneous Dirichlet boundary conditions while remaining smooth inside the domain.

\section{Error Evaluation Protocols}
\label{app:error_evaluation}

This appendix describes the mesh-based error metrics used to evaluate
the generalization performance of the proposed models.

Since all neural operator models produce predictions at discrete sensor locations,
both the predicted solution $\hat{u}$ and the reference solution $u$
are first mapped onto the computational domain $\Omega$
using the finite element infrastructure provided by DOLFINx.
Specifically, the point-wise values are interpolated onto
a piecewise-linear (P1) finite element function space defined on the underlying mesh.

The relative $L^2(\Omega)$ and $H^1(\Omega)$ errors are then computed via numerical
integration over the mesh:
\begin{equation}
\text{Rel. } L^2 \text{ error} =
\frac{\|u - \hat{u}\|_{L^2(\Omega)}}{\|u\|_{L^2(\Omega)}},
\qquad
\text{Rel. } H^1 \text{ error} =
\frac{\|u - \hat{u}\|_{H^1(\Omega)}}{\|u\|_{H^1(\Omega)}},
\end{equation}
with the norms defined as
\begin{equation}
\|v\|_{L^2(\Omega)}^2 = \int_{\Omega} |v(x)|^2 \, dx,
\qquad
\|v\|_{H^1(\Omega)}^2 =
\int_{\Omega} \left( |v(x)|^2 + |\nabla v(x)|^2 \right) dx.
\end{equation}
All integrals are evaluated using the assembled weak forms in DOLFINx,
ensuring consistent error measurement across different geometries.

\end{document}